%% file: main.tex
\documentclass[preprint,11pt]{elsarticle}

\usepackage{geometry, amsmath, dsfont}
\usepackage{graphicx} % Required for inserting images
\usepackage{xcolor}
\usepackage{url}
\usepackage{siunitx}
\usepackage{rotating}
\usepackage{booktabs} 
\usepackage{verbatim}
\usepackage{amssymb}
\usepackage{tikz}
\usetikzlibrary{calc}
\usepackage{float}
\usepackage{subcaption}
\usepackage{multirow} 
\usetikzlibrary{decorations.pathreplacing}
\usepackage{pgfplots}
\usepackage{graphicx}
\usepackage{gensymb}
\usepackage{makecell}
\usepackage{hyperref}

\usepackage{natbib}
 \bibpunct[, ]{(}{)}{,}{a}{}{,}%

\newcommand{\Nd}{\ensuremath{\mathds{N}}}
\newcommand{\Pd}{\ensuremath{\mathds{P}}}

\newcommand{\Ed}{\ensuremath{\mathds{E}}}
\newcommand{\Ac}{\ensuremath{\mathcal{A}}}

\newcommand{\Cc}{\ensuremath{\mathcal{C}}}
\newcommand{\Dc}{\ensuremath{\mathcal{D}}}
\newcommand{\Fc}{\ensuremath{\mathcal{F}}}

\newcommand{\Lc}{\ensuremath{\mathcal{L}}}

\newcommand{\Pc}{\ensuremath{\mathcal{P}}}

\newcommand{\Tc}{\ensuremath{\mathcal{T}}}

\newcommand{\Vc}{\ensuremath{\mathcal{V}}}

\newcommand{\Xc}{\ensuremath{\mathcal{X}}}

\input{new_command}

\DeclareMathOperator*{\argmin}{arg\,min}

\begin{document}

\begin{frontmatter}

%\RUNAUTHOR{Giliberto, Paradiso, and Wozabal}
%\RUNTITLE{Stochastic Programming for Dynamic Temperature Control of Refrigerated Road Transport}
\title{Stochastic Programming for Dynamic Temperature Control of Refrigerated Road Transport}

%\ARTICLEAUTHORS{%
% \AUTHOR{Francesco Giliberto, Rosario Paradiso, David Wozabal}
% \AFF{Vrije Universiteit Amsterdam, Amsterdam, Netherlands, \EMAIL{d.wozabal@vu.nl}, \EMAIL{d.wozabal@vu.nl}, \EMAIL{d.wozabal@vu.nl}}
%}% end of the block

\author[FG]{Francesco Giliberto}
\ead{f.giliberto@vu.nl}
\author[FG]{Rosario Paradiso}
\ead{r.paradiso@vu.nl}
\author[FG]{David Wozabal}
\ead{d.wozabal@vu.n}
\address[FG]{Department of Operations Analytics, Vrije Universiteit Amsterdam}

\begin{abstract}
     Temperature control in refrigerated delivery vehicles is critical for preserving product quality, yet existing approaches neglect critical operational uncertainties, such as stochastic door opening durations and heterogeneous initial product temperatures. We propose a framework to optimize cooling policies for refrigerated trucks on fixed routes by explicitly modeling these uncertainties while capturing all relevant thermodynamic interactions in the trailer. To this end, we integrate high-fidelity thermodynamic modeling with a multistage stochastic programming formulation and solve the resulting problem using stochastic dual dynamic programming. In cooperation with industry partners and based on real-world data, we set up computational experiments that demonstrate that our stochastic policy consistently outperforms the best deterministic benchmark by $35$\% on average while being computationally tractable. In a separate analysis, we show that by fixing the duration of temperature violations, our policy operates with up to $40$\% less fuel than deterministic policies. Our results demonstrate that pallet-level thermal status information is the single most crucial information in the problem and can be used to significantly reduce temperature violations. Knowledge of the timing and length of customer stops is the second most important factor and, together with detailed modeling of thermodynamic interactions, can be used to further significantly reduce violations. Our analysis of the optimal stochastic cooling policy reveals that preemptive cooling before a stop is the key element of an optimal policy. These findings highlight the value of sophisticated control strategies in maintaining the quality of perishable products while reducing the carbon footprint of the industry and improving operational efficiency.
\end{abstract}


\end{frontmatter}

\section{Introduction} \label{sec:intro}
% Cold chain market growth
Temperature-regulated supply chains, sometimes referred to as \emph{cold chains}, are essential in agricultural and pharmaceutical logistics. The global cold chain logistics market, valued at $313.3$ billion in 2024, is projected to grow to $410.7$ billion by 2028, driven by globalization and the growth of e-commerce \citep{statista, song2022research}.

% Efficient refrigeration is key
In cold chains, refrigeration and temperature monitoring are crucial to ensure quality and preserve the shelf life of products \citep{li2012dynamic}. Despite advances in equipment, such as refrigeration units and sensor technologies, the risk of product spoilage due to inadequate refrigeration remains significant \citep{ aiello2012simulation, konovalenko2021real}. This issue is particularly severe during the distribution phase, when trucks visit numerous customers within the same route \citep{AUNG2014198}, and operators frequently open the vehicle's doors, resulting in significant heat infiltration \citep{tassou2009food}. To counteract this and maintain more stable product temperatures, a refrigeration unit installed in the truck's trailer is typically activated when sensors detect a deviation from the admissible temperature range \citep{li2010optimal, HUANG2017576}.

While refrigeration is essential for the transport of perishable goods, it also represents a significant share of global energy consumption \citep{larsen2007potential}. According to \cite{deng2023pickup}, the fuel consumption of a refrigerated truck can be up to $1.5$ times higher than that of a non-refrigerated truck. As a result, cold chain operations, including transportation, contribute about $1$\% to global CO$_2$ emissions \citep{taher2021thermal}. These statistics highlight the necessity for innovative technical and planning solutions that promote efficient and sustainable transport operations \citep{hillier2015introduction,akkacs2022om}.

% A short statement of what we do
This paper examines a novel method to control refrigeration units in refrigerated road transport in order to keep the temperature of the products in a permissible range. The decisions are based on detailed modeling of the thermodynamic interactions, as well as data on the duration and timing of door openings, taking into account uncertainties in both of these aspects.

Challenges related to the distribution of perishable products have been investigated in various academic disciplines. A significant area of research in vehicle routing and related logistics problems focuses on extending existing models to incorporate additional requirements in cold chains. For example, \cite{rong2011optimization, DEKEIZER2017535} integrate a quality decay function in network design problems where the decision maker must decide on the production and distribution of perishable products and the temperature level at each storage facility. \cite{hsu2007vehicle} examine a vehicle routing problem with time windows, incorporating stochastic travel times and perishable products. Their study investigates how uncertain delivery times impact product quality.  
\cite{stellingwerf2018reducing, stellingwerf2021quality} incorporate a similar quality decay function in a routing problem where perishable products must be distributed from a storage facility to supermarkets, taking into account quality decay and energy consumption due to refrigeration and routing decisions. Although the above works embed quality decay functions that depend on temperature, they do not model the working regime of the refrigeration unit and assume that the trailer temperature will increase every time doors are opened and decrease during driving at a given fixed rate. 

\cite{deng2023pickup} focus on a deterministic pickup and delivery routing problem for temperature-sensitive products, considering incompatibility constraints between products with different properties and storage temperatures. They consider a quality-dependent selling price but simplify the degradation process by assuming linear perishability over time and also do not model the product's temperature evolution, assuming that the refrigeration unit maintains a constant temperature throughout the entire trip. Incompatibility constraints are also studied in works where vehicles are equipped with multiple compartments to transport products with conflicting temperature requirements \citep{CHEN20092311, OSTERMEIER2021799}. \cite{LIN2024, LIN2023103084} investigate the effect of precooling operations for post-harvest fruits and vegetables before the products enter the temperature controlled chain. The resulting problem is a deterministic heterogeneous fleet vehicle routing problem with time windows and temperature-dependent product quality decay. Their results underscore the importance of precooling in the early stages of the chain to preserve product quality.

Other studies explore more integrated decision-making problems within cold chains, combining production, inventory, and routing decisions. However, these studies feature less sophisticated models of perishability and quality. For example, \cite{CHEN20092311} study a production scheduling routing problem with uncertain demand. This problem involves determining when to produce and how to route vehicles, assuming that products have limited shelf-life and their value decays continuously over time. An essential limitation of this study is that the shelf-life is assumed to be constant and unaffected by external temperature or the number of stops. Similar simplifying assumptions regarding perishability and shelf-life are made in the production and inventory problems studied by \cite{crama2018, ALVAREZ2020511, ALVAREZ2022102667}, where the value of a product is age-dependent, and its decay is not necessarily linear.

All the aforementioned works either neglect cooling decisions or assume that the refrigeration unit can establish the desired temperature regardless of external conditions without modeling its working regimes in detail. Moreover, most papers consider deterministic problems or uncertainty limited to routing times. 

More accurate models to describe the working conditions of the refrigeration unit, its interactions with external and internal environments, and their energy consumption are studied in another stream of literature. These works consider models that use thermodynamic laws to describe temperature evolution in refrigerated environments. For example, \cite{de2015experimental} study heat infiltration in a trailer during door openings under various conditions, while \cite{hovgaard2011power, larsen2007potential} propose models to estimate the energy consumption needed to maintain the temperature in a refrigerated room. Other works use these models to improve truck designs for better insulation, e.g., \cite{tso2002experimental} analyze the impact of using different air curtain configurations, while \cite{zhang2018computational} focus on door sizes and locations. 

Other studies optimize the control schemes of the refrigeration unit \citep{li2010optimal,pei2021robust,song2022research,HUANG2017576,li2012dynamic}. However, these schemes are deterministic and agnostic to the routing of the trucks: they regulate the refrigeration unit based solely on the current air temperature without considering routing information, such as when the vehicle will stop with open doors and for how long.

Our work addresses these gaps in the literature by integrating thermodynamic modeling into a transportation context: we propose new dynamic cooling policies to control the refrigeration unit while considering the vehicles' routing and related uncertainties as input to make decisions. Our goal is to minimize the product's temperature deviations from allowed temperature ranges. To develop and validate our methodology, we collaborated with three companies active in the field of cold chain logistics: \emph{Mandersloot}, a third-party logistics provider specializing in refrigerated transportation; \emph{Tcomm Telematics}, which specializes in IT support and sensor technologies; and \emph{Ortec}, which specializes in decision support software. Overall, our work makes the following contributions to the current literature.
\begin{enumerate} 
    \item We propose the first stochastic cooling policy to control the refrigeration unit of a refrigerated vehicle, considering the uncertainties related to the duration of door opening and initial temperatures of the products loaded at each stop, which represent the main sources of uncertainty affecting loaded products' temperatures for a given route. The resulting problem is a large, multistage, stochastic programming problem, which we solve using stochastic dual dynamic programming (SDDP).
    
    \item We model the complete thermodynamic interactions in the trailer, striking a balance between physical accuracy and mathematical tractability. In particular, we do not only model the temperature of the air in the trailer but also explicitly account for the temperature of the products and the cooling fluid as well as two-way interactions between the air and every single product. This detailed modeling allows us to accurately track and consequently minimize excursion of product temperatures outside of allowed ranges.
    
    \item With the help of our industrial partners, we set up a detailed case study in which we investigate the impact of relevant parameters on a large number of problem instances. In particular, we vary route profiles, heat transfer coefficients of products, and the capacity of the cooling unit to create an ensemble of diverse, yet realistic settings that we use to benchmark our stochastic cooling policy out-of-sample against a deterministic lookahead policy, two simple myopic on/off heuristics from practice, and a perfect information benchmark.

    The results demonstrate that, across all instances, the heuristics perform significantly worse than the deterministic solution, which has access to detailed thermodynamic models, and that the stochastic solutions further outperform this latter policy by around $35\%$. We also show that the stochastic solution is only $7\%$ worse than the perfect information benchmark.
    
    \item By studying the results of our case study, we arrive at several structural insights and managerial implications:
    \begin{enumerate}
        \item Measurement of the temperature of the products, not only the air in the trailer, is the key to effective cooling decisions and represents the most important information in the problem. This speaks for the deployment of more sophisticated monitoring systems in trailers.

        \item Precooling before stops is important and the amount and timing of precooling are the main differences between our stochastic policy and other policies. This insight can be used as a starting point for the design of simple yet effective cooling policies that are based on product as well as air temperatures and explicitly manage a thermal buffer. Based on this insight, we demonstrate that such heuristics have the potential to outperform deterministic planning and come close to the optimum found by the stochastic programming solution.

        \item When imposing a limit on the maximal duration of temperature excursions, we show that our stochastic policy requires up to $40\%$ less fuel than the deterministic benchmark, demonstrating that there is scope for a significant reduction of the climate footprint of the industry, simply by optimized use of existing equipment.
        
        \item The information on the length of stops is the most impactful random factor in the problem and, in particular, is significantly more important than the knowledge of the initial temperatures of the products at the time of loading.

        \item Surprisingly, stronger cooling units do not lead to significantly fewer violations of temperature bounds and are only useful in routes with a large number of short-paced stops.
    \end{enumerate}
\end{enumerate}

The remainder of the paper is organized as follows: Section~\ref{sec:model} provides a detailed description of the problem under investigation and our thermodynamic modeling and casts the problem of taking optimal cooling decisions as a stochastic dynamic optimization problem. In Section~\ref{sec:policies}, we present four different cooling policies, drawing on relevant research from the stochastic programming literature. Finally, in Section~\ref{sec:case_study}, we present a realistic case study, including out-of-sample results obtained from a large set of instances. This section also offers practical recommendations and suggests potential directions for future research. Section \ref{sec:conclusion} concludes the paper.

{\bf Notation}: All random quantities are defined on a common probability space $\left(\Omega, \sigma, \Pd \right)$. For a random variable $X$, we write $X \triangleleft \Ac$ to indicate that $X$ is measurable with respect to the sigma algebra $\Ac$. Moreover, we define as $\Fc = (\Fc_1, \dots, \Fc_S)$, the natural filtration for the stochastic process $\xi = (\xi_1,, \dots, \xi_t)$, i.e., the smallest sigma algebras such that $\xi_t \triangleleft \Fc_t$. All equations and inequalities involving random quantities are meant to hold almost surely (a.s.) unless otherwise stated. For an integer $n \in \Nd$, we write $[n]$ for the set $\{1, \dots, n\}$.

%%%%%%%%%%%%%%%%%%%%%%%%%%%%%%%%%%%%%%%%%%%
%%% 2. Thermodynamic and Decision Model %%%
%%%%%%%%%%%%%%%%%%%%%%%%%%%%%%%%%%%%%%%%%%%
\section{Thermodynamic and Decision Model} \label{sec:model}
In this section, we present the decision model for controlling the refrigeration unit of a refrigerated truck. Section \ref{ssec:model_overview} provides an overview of our approach. Section \ref{ssec:thermodynamic_model} details the thermodynamic model that governs the relationship between the temperatures of different media inside the truck, together with the model for the unit’s power consumption. Finally, Section \ref{ssec:decision_model} introduces the stochastic programming formulation used to make cooling decisions.

%%% Section 2.1: Problem Description
\subsection{Problem Description} \label{ssec:model_overview}
We model the problem of controlling the operation of the refrigeration unit installed in a truck while picking up and delivering perishable products to a predefined set of customers on a predefined route as a stochastic programming problem. The goal consists of dynamically controlling the refrigeration unit, trying to keep the temperature of the products throughout the trip within their prescribed optimal temperature range. To this end, we assume that the duration of the stops at the customers' locations and the temperatures of the products at the time of loading are uncertain. These two factors are the most disruptive and have the highest potential to cause deviations from prescribed storage conditions and, thereby, product loss when there are deviations from nominal values. For the sake of simplicity, we assume that the travel times and weather conditions along the route are deterministic and known to the decision maker. However, also these factors could easily be modeled as stochastic by simple extensions of our model. 

\begin{figure}[t]
     \centering \includegraphics[width=1.0\textwidth]{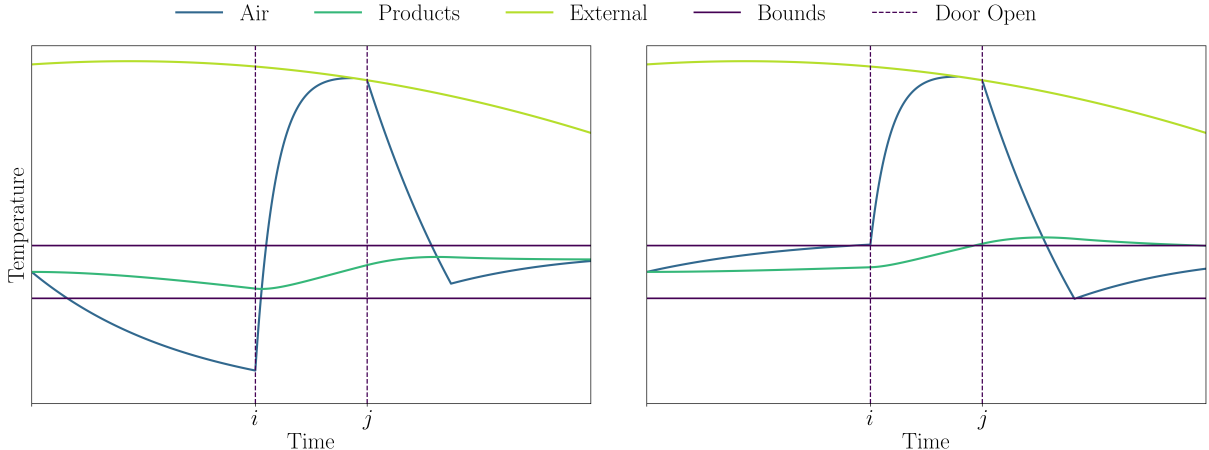}
    \caption{ {\small Temperature profiles of the air and products inside the trailer during a trip with a generic stop at a customer. The left panel shows a pre-emptive cooling strategy that cools the air and, consequently, the products in anticipation of the stop, while the right panel shows a purely reactive policy that starts cooling once the air temperature leaves a certain predefined range.}} 
    \label{fig:temperature_profile}
\end{figure}

Figure~\ref{fig:temperature_profile} illustrates the main thermodynamic effects governing the temperatures of the products and the air inside the trailer. In particular, the plot shows the temperature dynamics during a driving phase, a stop at a customer, and the subsequent driving phase under two different cooling strategies. When the truck starts its journey, both the products and the air are within the required temperature range (purple lines). In the left part of the figure, the refrigeration unit is immediately activated to preemptively cool before the doors open, gradually lowering the air temperature. Notice that the average temperature of the products also decreases, but at a slower rate due to their relatively higher heat capacity, making them harder to cool. In contrast, the right part of the figure illustrates a purely reactive policy, where cooling is only activated when the air temperature exceeds the allowed threshold. The latter, referred to as the 'on-off' rule, is typically adopted in practice \citep{li2010optimal, HUANG2017576}. 

When the truck stops at the customer at time $i$, the doors are opened, and the refrigeration unit cannot be used \citep{tso2002experimental, li2012dynamic}. This leads to a sharp increase in the temperature of the air and a gradual increase in the average temperature of the products. When the doors are closed at $j$, in both cases the refrigeration unit is activated to return the air temperature to an acceptable range. Although cooling the air inside the vehicle is relatively quick, cooling the products themselves takes longer. This delay is why a purely reactive policy fails to prevent temperature violations. In contrast, the preemptive policy can avoid this violation by building a thermal buffer in products that counteracts the effect of the influx of warm air.

The example illustrates the interplay between the air and products' temperature, but also the crucial role of cooling decisions in maintaining required storage conditions, especially during stops with uncertain durations. As we will show in Section \ref{sec:case_study}, preemptive cooling is instrumental in avoiding violations, and stochastic programming is able to find good preemptive cooling policies.

%%% 2.2 A Thermodynamic Model
\subsection{Thermodynamic Model}\label{ssec:thermodynamic_model}
We model the system's thermodynamics and the decision-making process over a discretized time horizon. The journey of the truck is represented as a sequence of points in time $\Tc$, in between which the truck is either in transit or stopped at a customer. We denote by $\Lc \subset \Tc$, points such that for $i \in \Lc$ the truck is waiting at a customer location with the doors open in the time interval $(i, i+1)$. Furthermore, let $\Dc \subset \Tc$ indicate the set of time points where the truck is driving, with $\Lc \cup \Dc = \Tc$. Finally, we denote by $\Delta_i$ the time elapsed between $i$ and $i+1$, expressed in seconds.

The truck transports a set of products whose composition changes over time due to handling operations at customers that include both loading and unloading of goods. At each time point, we model the temperatures of the air, the products, the cooling fluid that circulates within the refrigeration system, and the outside environment. To that end, we denote the set of products inside the truck at $i \in \Tc$ by $\mcl{P}_i$. At every $i \in \mcl{T}$, the temperature of the product $p \in \mcl{P}_i$ is denoted by $T_{pi}$ and should be maintained above a lower bound $\underline{T}_p$ and below an upper bound $\overline{T}_p$. \cite{goransson2018shelf} find that due to thermal inertia, the position of food products loaded on a pallet does not significantly affect the corresponding quality loss. Consequently, we assume homogeneous temperatures for products on the same pallet. Note that modeling the temperatures of the products requires information on the initial temperatures of the goods.

We assume that heat exchange occurs exclusively via the air in the truck, which interacts with the products, the external environment through the walls of the compartment and through the open doors during loading activities, as well as with the refrigeration system. We assume that the outside of the truck acts as a heat sink, i.e., that the truck does not influence the ambient temperature. 

The exchange of heat is governed by Newton's law of cooling \citep[see, e.g.,][]{lienhard2005heat}, which states that the thermal power $Q$ transferred between two media $a$ and $b$ is proportional to the difference in temperatures $T^a$ and $T^b$ of $a$ and $b$ as follows,
\begin{equation} \label{eq:newtons_law}
    Q = Ah (T^a - T^b),
\end{equation}
where $A$ is the area of contact [\unit{\meter\squared}] and $h$ [\SI{}{W\per\meter\squared\kelvin}] is the convective heat transfer coefficient that depends, among other things, on the properties of the media. Note that if $T^a < T^b$, then $Q<0$, which means that the thermal power is transmitted from $b$ to $a$. On the other hand, the heat dissipated or absorbed due to a variation of temperature of a medium $a$ between two time points $i$ and $i+1$ ($\Delta_i$ seconds apart) can be expressed as thermal power as follows,
\begin{equation} \label{eq:newtons_intertemporal}
    Q_i = C\frac{T_{i+1}^a - T_{i}^a}{\Delta_i},
\end{equation}
where $C$ is the heat capacity of the medium [\unit{\joule\per\kelvin}] \cite[e.g., ][]{zhang2018computational}. Note that \eqref{eq:newtons_intertemporal} is a linear discrete-time approximation of a continuous-time differential equation, which only works well if time intervals $\Delta_i$ are relatively short since it assumes that the rate of thermal power exchange remains constant during the interval.

We start by relating \eqref{eq:newtons_law} and \eqref{eq:newtons_intertemporal}, linking the change in air temperature between two time points $i$ and $i+1$ and the corresponding thermal power transfer obtaining,
\begin{equation}\label{eq:general_q}
    Q_i = C \frac{T^{\text{air}}_{i+1} - T^{\text{air}}_i}{\Delta _i},  \quad \forall i \in \mcl{T}\setminus \left\{|\mcl{T}|\right\},
\end{equation}
where $Q_i$ is the thermal power exchanged in $i$ [\unit{\watt}], $C$ is the heat capacity of the air inside the trailer [\unit{\joule\per\kelvin}], while $T^{\text{air}}_i$ represents its temperature at $i \in \Tc$ [\unit{\kelvin}].

Due to the principle of energy conservation, the thermal transfer $Q_i$ must equal the sum of the thermal powers exchanged in $(i, i+1)$ between the air and the external environment ($Q^{\text{ext-air}}_i$), between the air and each product ($Q^{\text{prod-air}}_{pi}$), and between the air and the refrigeration unit ($Q^{\text{cu-air}}_i$).

Consequently, we get
\begin{equation}\label{eq:energy_balance}
    Q_i = Q^{\text{ext-air}}_i + \sum_{p \in \mcl{P}_i}Q^{\text{prod-air}}_{pi} + Q^{\text{cu-air}}_i, \quad \forall i \in \mcl{T}.
\end{equation}

Due to \eqref{eq:newtons_law}, the thermal power exchanged with the external environment is
\begin{equation}\label{eq:q_air_ext}
    Q^{\text{ext-air}}_i  = \alpha_i\left(T^{\text{ext}}_i - T^{\text{air}}_i\right), \quad \forall i \in \mcl{T},
\end{equation}
where $T^{\text{ext}}_i$ is the temperature of the external environment at $i$ [\unit{\kelvin}]. The coefficient $\alpha_i$ is a constant that depends on the status of the doors at $i$. If the doors are closed, $\alpha_i$ is equal to the product between the wall's transmittance $U$ [\SI{}{W\per\meter\squared\kelvin}] and the exterior area of the trailer $A$ [\unit{\meter\squared}], where $U$ incorporates the contribution due to thermal conduction through the wall in addition to the phenomenon of convection with the external environment. When the doors are open, $\alpha_i =  UA' + vc'$,  where $A'$ is the surface of the trailer excluding the doors, $v$ is the mass of outside air entering the cargo [\unit{\kilogram\per\second}], and $c'$ specific heat of the internal air [\SI{}{J\per\kilogram\kelvin}]
\citep{zhang2018computational}.

Similarly, the thermal power exchanged between the air and each product is
\begin{equation}\label{eq:q_air_prod}
    Q^{\text{prod-air}}_{pi} =  \beta_{pi}\left( T^{\text{}}_{pi} - T^{\text{air}}_i \right), \quad  \forall i \in \mcl{T}, \quad p \in \mcl{P}_i,
\end{equation}
where $\beta_{pi} = h_p A_p$, with $h_p$ representing the heat transfer coefficient of product $p$ and $A_p$ its surface area in contact with the air. Finally, the thermal power exchanged between the air and the refrigeration unit is given by 
\begin{equation}\label{eq:q_air_cu}
    Q^{\text{cu-air}}_i = \gamma_i\left(\tcu_i - \ttruck_i\right), \quad \forall i \in \mcl{T},
\end{equation}
where $\tcu_i$ is the temperature of the refrigerant fluid at time point $i$ and $\gamma_i$ is the overall heat transfer coefficient between evaporator and the air which can be assumed to be constant if the evaporator fan is working at a constant speed \citep{larsen2007potential}. This parameter is set equal to $0$ for each $i \in \Lc$, since the unit is off. By combining \eqref{eq:general_q}~--~\eqref{eq:q_air_cu}, we obtain the following equation for the evolution of air temperature over time
\begin{equation}\label{eq:air_temperature_evolution}
	\ttruck_{i+1} = \ttruck_{i} + \frac{\Delta_i}{C}\left(\alpha_i\left(\tamb_i - \ttruck_i\right) 
	+   \sum_{p \in \mcl{P}_i}\beta_{pi}\left(\tprod_{pi} -\ttruck_i\right) + \gamma_i\left(\tcu_i - \ttruck_i\right)\right), \quad \forall i \in \mcl{T} \setminus \{|\mcl{T}|\}.
\end{equation}

Correspondingly, the thermal power exchange in $(i, i+1)$ between the air and a product $p \in \mcl{P}$ can be related to the evolution of each product temperature in the following way 
\begin{equation}\label{eq:product_energy_balance}
    Q^{\text{prod-air}}_{pi} =   C_p\frac{T_{p,i+1} - T_{pi}}{\Delta_i},    \quad \forall i \in \mcl{T}\setminus \{|\mcl{T}|\}, \quad p \in \mcl{P}_i \cap \mcl{P}_{i+1},
\end{equation}
where $C_p$ is the heat capacity of the product $p$. Note that by restricting the products to $\mcl{P}_i\cap\mcl{P}_{i+1}$, we assume that if a product is unloaded in $(i,i+1)$, it is unloaded at the start of the time interval. By combining equations \eqref{eq:q_air_prod} and \eqref{eq:product_energy_balance}, we obtain the evolution of the temperature of each product over time as a function of the temperature of the air
\begin{equation}\label{eq:product_temperature_evolution}
    T_{p,i+1}  =T_{pi} + \frac{\Delta_i\beta_{pi}}{C_p}\left(T^{\text{air}}_i - T^{\text{}}_{pi}\right),   \quad \forall i \in \mcl{T}\setminus \{|\mcl{T}|\}, \quad p \in \mcl{P}_i \cap \mcl{P}_{i+1}.
\end{equation}

To control the internal conditions of the truck, we dynamically regulate the temperature of the refrigerant fluid at the evaporator of the refrigeration system $(\tcu)$, which is in direct contact with the air \eqref{eq:q_air_cu}. To model the operational limits of the unit, we ensure that the consumption of electrical power $W_i$ resulting from the model's decisions does not exceed the actual capacity $\overline{W}$ [\unit{\watt}], i.e.,
\begin{align}\label{eq:watt_ub}
    W_i \in [0,\overline{W}], \quad \forall i \in \Dc.
\end{align}

The consumption $W_i$ needed to drive the process depends both on the amount of heat absorbed by the refrigerant \eqref{eq:q_air_cu} and on the unit's coefficient of performance (COP), which in turn depends on the difference between the internal and external air temperature. \cite{hovgaard2011power, larsen2007potential} model the ensuing complex interactions considering the physical properties of the evaporator by
\begin{align}\label{eq:cooling_power}
    W_i = \theta_{i1} {\tcu_i}^2 + \theta_{i2} \tcu_i + \theta_{i3} \tcu_i\ttruck_i + \theta_{i4} T^{\text{air}}_i, \quad \forall i \in \Dc, 
\end{align}
where the coefficients $\theta_{ij}$ model the thermodynamic properties of the system, and the details of their calculation are provided in the Appendix \ref{app:cooling_unit}. It is reasonable to assume that once $(\tcu)$ is set, the system is able to reach it immediately, regardless of the previous value \citep{larsen2007potential}. Therefore, in \eqref{eq:cooling_power}, the instantaneous power to be supplied to the unit does not depend on the state of the system at the previous instant, but only on the current state.

Unfortunately, \eqref{eq:cooling_power} is not a convex function in our decision variables $\ttruck_i$ and $\tcu_i$. We therefore construct a convex approximation as a maximum of $K$ affine functions using the least squares partition algorithm proposed in \cite{magnani2009convex}. More specifically, we estimate a set of parameters ($\phi \in \mathbb{R}^{3 \times |\Dc| \times K}$) for each possible outdoor temperature value recorded during the route using a grid of values for $(\ttruck_i, \tcu_i)$ as well as the corresponding values for $W_i$ according to \eqref{eq:cooling_power}. The power consumption constraint can then be rewritten as
\begin{align}\label{eq:cooling_power_approx}
    W_i &\geq \phi_{0, i, k} + \phi_{1, i, k} \, \tcu_i + \phi_{2, i, k} \, \ttruck_i  \quad \forall k \in \{1, \ldots, K\}, \quad \forall i \in \Dc,
\end{align}
where \(k\) denotes the index of the linear function corresponding to the external temperature at time \(i\) with parameters \(\phi_{0, i, k}, \phi_{1, i, k}, \phi_{2, i, k}\). Depending on the properties of the refrigerant and the system’s technical specifications, it is reasonable to impose a lower bound on the variable $\tcu$
\begin{align}\label{eq:cooling_limits}
    \tcu_i \geq \Gamma, \quad \forall i \in \Dc.
\end{align}

Finally, we add a constraint to prevent heating, i.e., that the refrigerant temperature exceeds the internal air temperature \citep{larsen2007potential}
\begin{align}\label{eq:only_cooling}
    \tcu_i \leq \ttruck_i, \quad \forall i \in \Dc.
\end{align}
%Note that with an appropriate energy consumption function for a heat pump, our model could be generalized to incorporate heating as well. However, here we focus exclusively on cooling, which is in line with most real-world settings.
While our model could potentially be extended to include heating by incorporating an appropriate energy consumption function for the heating regime, our current focus remains exclusively on cooling. %This choice is in line with the majority of real-world applications.

%%% 2.3 The stochastic dynamic model
\subsection{The Stochastic Dynamic Model} \label{ssec:decision_model}
In this section, we describe the multi-stage stochastic programming model used to control the refrigeration unit installed in the truck. Planning is carried out per route, with a loaded truck departing from the depot at time $0$. We do not model a loading phase at the depot, and the optimization ends at the last stop $S$, where the truck is fully emptied. We assume that the route of the truck is planned offline and that the travel times, the cargo composition $\Pc_i$, as well as the outside temperature $\tamb_i$ are deterministic and known throughout the planning horizon. All of these assumptions could, in principle, be relaxed at the cost of more involved models. 

At each stop, the truck is partially unloaded and potentially loaded with new products before leaving for the next stop on the route, resulting in a pattern of alternating cargo handling and driving phases. We assume that the trucking company is responsible for the temperature of the goods until they arrive at their destination, but handling operations are the responsibility of the customers receiving the products. This, in particular, implies that we do not keep track of any violations that may occur during the unloading of the products, and as a result, we do not model the unloading process at the final stop. Furthermore, we update $\mcl{P}_i$ so that loaded products enter at the end of handling, while unloaded products are removed at the beginning.

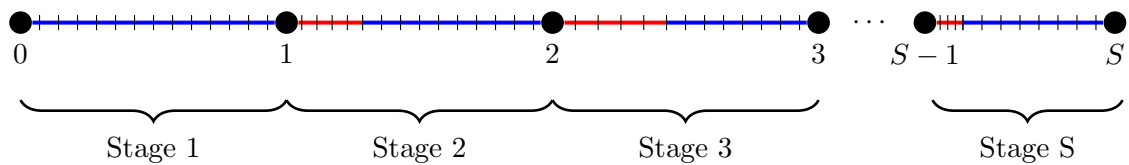
\begin{figure}[t]
    \begin{center}
        \begin{tikzpicture}[node distance=1.5cm, auto]
            % Nodes - only large dots at stage boundaries
            \node (1) [circle, fill=black, inner sep=3pt, label=below:0] {};
            \node (2) [right of=1, xshift=2cm, circle, fill=black, inner sep=3pt, label=below:$1$] {};
            \node (3) [right of=2, xshift=2cm, circle, fill=black, inner sep=3pt, label=below:$2$] {};
            \node (4) [right of=3, xshift=2cm, circle, fill=black, inner sep=3pt, label=below:$3$] {};
            \node (5) [right of=4, xshift=-0.8cm] {$\cdots$};
            \node (6) [right of=5, xshift=-0.8cm, circle, fill=black, inner sep=3pt, label=below:$S-1$] {};
            \node (7) [right of=6, xshift=1cm, circle, fill=black, inner sep=3pt, label=below:$S$] {};
            
            % Lines between dots - split into blue and red segments
            \draw[blue, line width=1.5pt] (1) -- (2);
            \foreach \x in {0.25,0.5,...,3.5} % Every 5mm along blue line
                \draw[thin] ($(1) + (\x,-0.1)$) -- ($(1) + (\x,0.1)$);
            
            \draw[red, line width=1.5pt] (2) -- ++(1cm,0);
            \foreach \x in {0.2,0.4,...,1} % 10 equally spaced lines
                \draw[thin] ($(2) + (\x,-0.1)$) -- ($(2) + (\x,0) + (0,0.1)$);
            
            \draw[blue, line width=1.5pt] ($(2)+(1cm,0)$) -- (3);
            \foreach \x in {1,1.25,...,3.5} 
                \draw[thin] ($(2)+(\x, -0.1)$) -- ($(2)+(\x, 0.1)$);
            
            \draw[red, line width=1.5pt] (3) -- ++(1.5cm,0);
            \foreach \x in {0.3,0.6,...,1.5} 
                \draw[thin] ($(3) + (\x,-0.1)$) -- ($(3) + (\x,0.1)$);
            
            \draw[blue, line width=1.5pt] ($(3)+(1.5cm,0)$) -- (4);
            \foreach \x in {1.5,1.75,...,3.5} 
                \draw[thin] ($(3) + (\x,-0.1)$) -- ($(3) + (\x,0.1)$);
            
            \draw[red, line width=1.5pt] (6) -- ++(0.5cm,0);
            \foreach \x in {0.1,0.2,...,0.5} 
                \draw[thin] ($(6) + (\x,-0.1)$) -- ($(6) + (\x,0.1)$);
            
            \draw[blue, line width=1.5pt] ($(6)+(0.5cm,0)$) -- (7);
            \foreach \x in {0.5,0.75,...,2.25} 
                \draw[thin] ($(6)+(\x,0.1)$) -- ($(6)+(\x,-0.1)$);
            
            % Braces
            \draw [decorate,decoration={brace,amplitude=8pt,mirror,raise=15pt},line width=1pt](0,-0.5) -- (3.5,-0.5) node [black,midway,yshift=-1.5cm] {Stage 1};
            \draw [decorate,decoration={brace,amplitude=8pt,mirror,raise=15pt},line width=1pt](3.5,-0.5) -- (7,-0.5) node [black,midway,yshift=-1.5cm] {Stage 2};
            \draw [decorate,decoration={brace,amplitude=8pt,mirror,raise=15pt},line width=1pt](7,-0.5) -- (10.5,-0.5) node [black,midway,yshift=-1.5cm] {Stage 3};
            \draw [decorate,decoration={brace,amplitude=8pt,mirror,raise=15pt},line width=1pt](12,-0.5) -- (14.5,-0.5) node [black,midway,yshift=-1.5cm] {Stage S};
        \end{tikzpicture}
    \end{center}
    \caption{\label{fig:time_structure} {\small The time structure of the multi-stage stochastic program: the truck starts at the depot at time $0$ and then loads/unloads refrigerated goods at $S$ stops. The time between the stops are the stages of the stochastic program, each consisting of a driving phase and all but the first of a loading/unloading phase. Within each stage, there is a finer temporal structure used to model thermodynamic interactions and decisions about cooling.}}
\end{figure}

We continue our discussion with a description of the temporal structure of the dynamic program depicted in Figure \ref{fig:time_structure}. As described in the previous section, we use our thermodynamic model to track the temperatures of air, goods, and cooling fluid at all time points $i \in \Tc$. Additionally, we impose a coarser time structure that nests the fine structure $\Tc$ and models the information flow in the stochastic dynamic program. In particular, we define $S \in \Nd$ decision stages and assign each $i$ in $\Tc$ to one stage, denoting the set of $i$ that belongs to stage $s$ by $\Tc_s$ so that the decision stage $s$ contains all time points $\Tc_s\subseteq \Tc$ between arriving at stop $s$ and arriving at stop $s+1$. In line with the discussion above, we partition $\Tc_s=\Lc_s \cup \Dc_s$ into $i \in \Lc_s$ which models the cargo handling phase and $i\in \Dc_s$ which models the truck while driving. Furthermore, we define $i_0(s)$ as the first $i\in \Tc_s$ in the driving phase of stage $s$. %and $\Ic_s$ as the inventory of products that are transported in period $s$.

Randomness enters our model in the form of random door opening times $O_s$ and initial temperatures $T_{pi_0(s)}$ of all newly loaded products $p \in \Pc_{i_0(s)}\setminus \Pc_{i_0(s-1)}$ at each stop $s$, excluding the first one. These two parameters have the highest potential to cause significant deviations from planned cooling schedules and are revealed at the beginning of each stage, i.e., when the cargo handling phase starts. As a result, our random data in stage $s$ can be described as $\xi_s = ((T_{pi_0(s)})_{p: p\in \Pc_{i_0(s)}\setminus \Pc_{i_0(s-1)}}, O_s)$, which realizes at stage $s$. Note that the dimension of $\xi_s$ varies from stage to stage depending on how many new products are loaded. We denote by $\Fc=(\Fc_1, \dots, \Fc_S)$ the $\sigma$-algebra generated by this process. 

The random door opening times $O_s$ have implications for the discretization of the time of cargo handling periods. Although time intervals during driving $\Delta_i$ with $i \in \Dc_s$ can be chosen to be of fixed length $\Delta_D$ for deterministic driving times, the same is not possible for the cargo handling periods. For this reason, to model $\Delta_i$ for $i \in \Lc_s$, we fix a number $|\Lc_s|$ of time periods upfront and independent of the realization of the random variables $O_s$ and then define $\Delta_i = O_s/|\Lc_s|$ for all $i \in \Lc_s$. This implies that the duration of each time step $\Delta_i$ while the doors are open is random. Note that this is the only parameter directly affected by random door opening times.

Since new information on $\xi_s$ only arrives when the truck stops to load and unload, the cooling decisions in all time instants during the drive between stops $s$ and $s+1$ are measurable with respect to $\Fc_s$ and, therefore, are taken at stage $s$. Hence, we have a time structure in stages that correspond to client stops and, on top of that, a fine grained time structure which enables frequent adjustments of the refrigeration unit during travel between stops.
    
The objective of controlling the refrigeration unit of the truck is to avoid that the temperature of the transported goods leave the prescribed storage conditions. In order to represent this goal as the objective function, we minimize the expected violation of the temperature targets. More specifically, we measure the violations by 
\begin{align}
    \nu_{pi} = \max( T_{pi} - \overline T_p, 0) + \max( \underline T_p - T_{pi}, 0)
\end{align}
for all instants $i \in \Tc_s$ and all products $p \in \mcl{P}_i$. By adding up positive and negative excursions from the permitted temperature range, we implicitly assume that the costs of violations are symmetric. Obviously, this assumption could be relaxed by assigning different weights to the two types of violation and also differentiating between products. Finally, in addition to the absolute temperature deviation, it is also important to consider the duration of this deviation \citep{rong2011optimization}, which is represented by $\Delta_i$. 

Our objective function sums up these cost terms for one scenario over all stages, all sub-stage time points, and products, resulting in
\begin{equation} \label{eq:objective}
    \Cc = \sum_{s \in [S]} \underbrace{\sum_{p \in \Pc_{i_0(s)}} \left( \sum_{i \in \Lc_s}  \Delta_i\nu_{pi} + \sum_{i \in \Dc_s} \Delta_D \nu_{pi} \right)}_{\Cc_s},
\end{equation}
where $\Cc$ is the function for a specific scenario and $\Cc_s$ is the cost in stage $s$ [\unit{\kelvin\cdot\second}].

Putting everything together and taking the expectation of all violation costs in the objective, we can therefore write the extensive form of our linear multi-stage stochastic optimization program as
\begin{subequations} \label{eq:extensive_form}
    \begin{align} 
      \min & \quad  \mathbb{E} \left[ \sum_{s \in [S]} \sum_{p \in \Pc_{i_0(s)}} \left( \sum_{i \in \Lc_s}  \Delta_i\nu_{pi} + \sum_{i \in \Dc_s} \Delta_D \nu_{pi} \right) \right] \nonumber \\
      \text{s.t.} & \quad \eqref{eq:air_temperature_evolution}, \eqref{eq:product_temperature_evolution}, & \text{[thermodynamic constraints]} \label{eq:thermo_dynamic_constraints}
      \\      
      & \quad \eqref{eq:watt_ub}, \eqref{eq:cooling_power_approx},
      \eqref{eq:cooling_limits}, & \text{[refrigeration unit's constraints]}
      \label{eq:unit_constraints}
      \\
      & \quad \eqref{eq:only_cooling}, & \text{[no heating]}
      \label{eq:no_heating}
      \\
      & \quad \ttruck_i, \tcu_i \triangleleft \Fc_{s},  & \forall i \in \Tc_s, \; \forall s \in [S] \label{eq:non_anticipativity_1}
      \\
      & \quad W_i \triangleleft \Fc_{s},  & \forall i \in \mathcal{D}_s, \; \forall s \in [S]\label{eq:non_anticipativity_2}
      \\
      & \quad T_{pi} \triangleleft \Fc_{s},  & \forall p \in \Pc_{i_0(s)}, \; \forall i \in \Tc_s, \; \forall s \in [S] \label{eq:non_anticipativity_3}
    \end{align}
\end{subequations}
where \eqref{eq:non_anticipativity_1}, \eqref{eq:non_anticipativity_2}, and \eqref{eq:non_anticipativity_3} enforce that all the decisions are non-anticipative. All constraints of the above problem hold almost surely, and the expectation is with respect to the image measure of $\xi$. 

In order to assess the differences in fuel consumption between the different policies, we additionally use a variant of model \eqref{eq:extensive_form}, which includes an additional constraint on the overall fuel budget $B$ [\unit{\liter}] to be used for powering the refrigeration unit 
\begin{align}\label{eq:fuel_budget_constraint}
    \sum_{s \in [S]} \sum_{i \in D_s} \frac{W_i\Delta_i}{\sigma} \leq B + \mu,
\end{align}
where $\sigma$ converts the electrical energy consumption of the refrigeration unit to induced fuel consumption \citep{stellingwerf2018reducing}. The nonnegative slack variable $\mu$ on the right-hand side is included with a large positive coefficient in the objective modeling penalty for not adhering to the fuel constraint.

Note that if we omit $\mu$ in \eqref{eq:fuel_budget_constraint}, the problem might not have relatively complete recourse. In fact, if the refrigeration system consumes all the available fuel budget until period $i$, the piecewise linear approximation \eqref{eq:cooling_power_approx} does not guarantee $W_i=0$ exactly if, even $\tcu_i=\ttruck_i$. This would therefore result in an infeasible problem in stage $s$ with $i \in \Dc_s$ given feasible decisions from stage $0$ to stage $s-1$. The above formulation ensures that the problem has relative recourse and effectively enforces the fuel constraint for the optimal policy. This is required since we use SDDP to solve the problem and want to avoid feasibility cuts. 

%%%%%%%%%%%%%%%%%%%
%%% 3. Policies %%%
%%%%%%%%%%%%%%%%%%%
\section{Cooling Policies} \label{sec:policies}
In this section, we outline the four cooling policies that we consider to address the problem defined in the previous section: a stochastic programming policy (SP) that directly solves problem \eqref{eq:extensive_form}, a rolling lookahead policy (RLP) that uses repeated deterministic planning based on a forecast of randomness $\xi$, but still using the full model of thermodynamic interactions, and lastly, two myopic 'on-off' policies that do not anticipate upcoming stops at all. The first myopic policy (H1) is unaware of the thermal state of the cargo and reacts only to the temperature of the air inside the trailer, while the second myopic policy (H2) also has information about each of the transported products.

\subsection{A Stochastic Programming Solution Based on SDDP}
The problem described in the last section is a multistage stochastic linear program whose size grows in the number of products and customer stops. We propose solving this problem using SDDP \citep{pintoPereira_sddp,philpott:2011,MinnerLohndorfWozabal2013,rebennack2016,downward2020,LoWo21}. Recently, SDDP has received a lot of attention due to its ability to solve stochastic optimization problems with many stages while avoiding the exponential growth of complexity that is traditionally associated with an increase in decision stages \citep[e.g.,][]{lan2020complexity}. 

We approximate the stochastic process $\xi=(\xi_1, \dots, \xi_S)$ governing temperatures of newly loaded products and door opening durations by a finite scenario lattice $\tilde \xi$ \citep{MinnerLohndorfWozabal2013, LoWo21}. This approach allows the stochastic program to be solved efficiently, resulting in a policy that can be applied out-of-sample, i.e., for scenarios not represented by the scenario lattices as is argued below.

Our approach solves \eqref{eq:extensive_form} as a stochastic dynamic program governed by its dynamic programming equations
\begin{align} \label{eq:dynamic_programming_equations}
    V_s(x_{s-1}, \xi_s) = \left\{ 
    \begin{array}{ll} 
       \min & \sum_{p \in \Pc_{i_0(s)}} \left( \sum_{i \in \Lc_s}  \Delta_i\nu_{pi} + \sum_{i \in \Dc_s} \Delta_D \nu_{pi} \right) + \Vc_s(x_s, \xi_s) \\
      \text{s.t.} & x_s \in \Xc_s(x_{s-1}, \xi_s)
    \end{array}
    \right.
\end{align}
where $x_s$ are the decisions that are made in stage $s$ that are relevant for stage $s+1$, i.e., the temperatures of the air and the products in all time periods associated with stage $s$ and
$$ \Vc_s(x_s, \xi_s) = \Ed\left[ V_{s+1}(x_s, \xi_{s+1}) | \xi_s \right]$$
is the post-decision value function \citep{powell:2011}.

The idea of SDDP is to iteratively refine approximations of the post-decision value functions. The resulting approximations $\tilde \Vc_s$ implicitly define a policy on the states that are part of the lattice. However, the policy can be lifted to arbitrary realizations of $\xi$. For that purpose, we define the rounded state as
$$ \tilde \xi_s(\xi_s) \in \argmin_{y \in \operatorname{supp}(\tilde \xi_s)} ||y - \xi_s||.$$

For a given state $\tilde x_{s-1}$ and a realization $\bar \xi_s$ of $\xi_s$, we then can make a decision $\tilde x_s$ by solving
\begin{align} \label{eq:apply_policy_oos}
    \begin{array}{ll} 
       \min & \sum_{p \in \Pc_{i_0(s)}} \left( \sum_{i \in \Lc_s}  \Delta_i\nu_{pi} + \sum_{i \in \Dc_s} \Delta_D \nu_{pi} \right) + \tilde \Vc_s(x_s, \tilde \xi_s(\bar \xi_s)) \\
      \text{s.t.} & \tilde x_s \in \Xc_s(\tilde x_{s-1}, \bar \xi_s).
    \end{array}
\end{align}
Note that, to avoid infeasibilities, we round from $\bar \xi_s$ to $\tilde \xi_s(\bar \xi_s)$ only when retrieving the value function but not when constructing the feasible set.

To estimate the objective of the policy, we then simulate a set of $N$ trajectories $(\bar \xi_s^n)_{s\in [S], n \in [N]}$ of $\xi$. For each of these trajectories, we set $\tilde x_0^n=x_0$ as the deterministic starting state of the problem and then recursively define $\tilde x_s^n$ using $\tilde x_{s-1}^n$, $\bar \xi_s^n$ and the problem \eqref{eq:apply_policy_oos} to find $\tilde x_s^n$. Finally, we calculate
$$ \tilde \Cc^{\text{SP}} = N^{-1} \sum_{n\in [N]} \sum_{s \in [S]} \tilde \Cc_s^n(\tilde x_s^n)$$
as the average cost of the policy on the $N$ samples, where $\Cc_s^n(\tilde x_s^n)$ is the cost in scenario $n$ and stage $s$ when taking the decision $\tilde x_s^n$. In this way, we evaluate the policy found by SDDP and encapsulated in $\tilde V_s$ \emph{out-of-sample} \citep[see][for further details]{LoWo21}.

\subsection{A Rolling Lookahead Policy}
We employ a rolling lookahead policy (RLP) as deterministic benchmark. The idea is to ignore the stochastic part of the optimization and solve the dynamic part of the problem by merely using the expectation of the stochastic process. This policy offers a simple solution heuristic for any dynamic stochastic optimization problem, which offers a powerful benchmark that often produces excellent results due to its ability to rely on forecasts of future realizations of $\xi$ and the fact that it is applied repeatedly, which gives RLP the opportunity to revise suboptimal decisions by adapting to new information. 

More specifically, RLP is a receding-horizon lookahead strategy which simulates the solution of a deterministic optimization over a rolling horizon, where future states of the stochastic process are replaced with their expectations. To calculate the value of the RLP, we again simulate $N$ trajectories $(\bar \xi_s^n)_{s\in [S], n \in [N]}$ of $\xi$ and set $\tilde x_0^n=x_0$ and $\tilde \xi_0^n=\xi_0$ the deterministic starting state for each trajectory $n \in [N]$. To obtain $\tilde x_{s'}^n$ from $\tilde x_{s'-1}^n$, we solve the following deterministic problem
\begin{align} \label{eq:RLP_policy}
    \begin{array}{ll} 
       \min & \sum_{s=s'}^S \sum_{p \in \Pc_{i_0(s)}} \left( \sum_{i \in \Lc_s}  \Delta_i\nu_{pi} + \sum_{i \in \Dc_s} \Delta_D \nu_{pi} \right) \\
      \text{s.t.} & (x_{s'}, \ldots, x_S) \in \Xc_{s'}(\tilde x_{s'-1}, \tilde \xi_{s'}) \times \cdots \times \Xc_{S}(x_{S-1}, \Ed[\bar \xi_S| \xi_{s'} = \tilde \xi_{s'}^n])
    \end{array}
\end{align}
which is simply \eqref{eq:apply_policy_oos} with the value function replaced by a lookahead to future costs based on the forecasts of $\xi_s$ given $\tilde \xi_{s}^n$. We then set $\tilde x_{s'}^n=x_{s'}^*$ for $x_{s'}^*$ the optimal solution for stage $s'$ in \eqref{eq:apply_policy_oos}. Subsequently, we resolve the problem for $(\tilde x_{s'}^n, \tilde \xi_{s'+1}^n)$ to obtain $\tilde x_{s'+1}^n$. In this way, we solve the problem \eqref{eq:RLP_policy} $N$ times to get a trajectory of solutions for $\tilde \xi ^n$. Each solution in a particular stage is based on the forecasts for the rest of the horizon, which are updated in a rolling fashion as new information becomes available.

The out-of-sample value of the policy is therefore defined as
$$ \tilde \Cc^{\text{RLP}} = N^{-1} \sum_{n\in [N]}\sum_{s \in [S]} \tilde \Cc_s^n(\tilde x_s^n)$$
with the $\tilde x_s^n$ generated as described above.

In summary, the RLP policy still features the full thermodynamic model and anticipates stops. However, RLP does not appreciate the stochastic nature of the problem and, therefore, the difference $\Cc^{\text{RLP}}-\Cc^{\text{SP}}$ can be interpreted as the value of the stochastic solution.

\subsection{Two Simple Heuristics}
To evaluate our proposed stochastic programming solution, we also benchmark against two myopic heuristics, i.e., policies that do not anticipate future events and therefore are unaware of upcoming handling operations and the resulting door openings. In both policies, we emulate the on-off feedback controller with hysteresis \citep{li2010optimal}, which is a commonly applied approach in practice due to its simplicity. Operationally, it consists of maintaining the internal temperature within a specific range by turning the system on or off. Specifically, this strategy involves operating the refrigeration unit at full capacity when the internal temperature reaches an upper allowable limit and turning it off when it reaches a lower limit \citep{HUANG2017576}. Both policies lack forward-looking decision-making but differ in their degree of awareness of temperature in the trailer: our first heuristic, H1, measures only the air temperature, whereas the second heuristic, H2, considers the entire thermal state, including tracking the temperature of each product.  

Among our policies, H1 is closest to current industry practice. We define its on/off state as
$$O_i^{\text{H1}} =
\begin{cases} 
1, & \text{if } \ttruck_i > \min_{p \in \Pc_{i}} \underline{T}_p + \lambda_1 (\overline{T}_p - \underline{T}_p), \\
1, & \text{if } O_{i-1}^{\text{H1}} = 1 \land  \ttruck_i > \max_{p \in \Pc_i} \underline{T}_p, \\
0, & \text{otherwise},
\end{cases}$$
where the parameter $\lambda_1$ controls how close the air temperature can reach the maximum allowed temperature for the goods before the refrigeration unit is turned on. Once the unit is on, it cools the air to the maximum of the minimum desirable temperatures of the products. In this way, the ideal storage condition is preserved if the unit operates in stable conditions without any door openings or other sudden temperature changes.

Our second myopic policy, H2, adds an additional control layer on top of H1, activating the system if one of the products approaches a critical state (sufficiently close to its maximum allowable limit) and deactivating it as soon as a product reaches its lower limit
$$O_i^{\text{H2}} =
\begin{cases} 
1, & \text{if } \Big( \ttruck_i > \min_{p \in \Pc_i} \underline{T}_p + \lambda_1 (\overline{T}_p - \underline{T}_p)  \lor   \max_{p \in \Pc_i} (T_{pi} - \underline{T}_p) - \lambda_2 (\overline{T}_p - \underline{T}_p) > 0\Big)  \\
   & \quad \land  \left( \min_{p \in \Pc_i} T_{pi} - \underline{T}_p > 0\right) , \\
1, & \text{if } O_{i-1}^{\text{H2}} = 1 \land  \min_{p \in \Pc_i} T_{pi} - \underline{T}_p > 0 \land  \ttruck_i > \max_{p \in \Pc_i} \underline{T}_p, \\
0, & \text{otherwise}.
\end{cases}$$
Hence, H2 activates the unit if the air overheats excessively (controlled by $\lambda_1$), or if one of the products exceeds the established thermal tolerance (controlled by $\lambda_2$). In both cases, the refrigeration unit is not switched on if one of the products reaches its lower thermal limit. Once the unit is turned on, it cools until the air temperature is below the maximal temperature limit of the products, and one of the products reaches its lower temperature limit. As we shall see in the next section, the additional consideration of product temperatures leads to clear benefits.

In both policies, if the unit is turned on, it operates at maximum capacity, that is,
\[
    W_i = O_i \cdot W_i^{\text{max}},
\]
with $W_i^{\text{max}}$ defined as
\[
    W_i^{\text{max}} = \min \left( \max_{k} \left\{ \phi_{0, i, k} + \phi_{1, i, k} \cdot \Gamma + \phi_{2, i, k} \cdot \ttruck_i \right\}, \, \overline{W} \right), \quad \forall i \in \mathcal{D}_s, \; \forall s \in [S]
\]
and the temperature of the cooling fluid $\tcu_i$ equals $\Gamma$ in the first case and is determined so that the power does not exceed $\overline{W}$ in the second case.

%%%%%%%%%%%%%%%%%%%%%
%%% 4. Case Study %%%
%%%%%%%%%%%%%%%%%%%%%
\section{Case Study}\label{sec:case_study}
In this section, we analyze the performance of the proposed policies using a real-life case study. Section \ref{ssec:instances_generation} discusses the employed data and the methods adopted to generate the problem instances. Section \ref{ssec:samples_generation} details the process of producing scenarios for door opening times and initial product temperatures. Finally, Section~\ref{ssec:numerical_results} presents a comprehensive analysis of the performance of different policies and the resulting cooling decisions, which allows for the derivation of structural insights and managerial implications.

%%% Section 4.1: Test Instance Generation
\subsection{Test Instance Generation}\label{ssec:instances_generation} 
Generating instances for our problem requires a combination of different types of information from different sources. These data include vehicle and refrigeration unit characteristics, route details, weather conditions, and specific thermodynamic parameters needed for the physical modeling of the system. 

We obtain this information through our cooperation with a consortium of industry partners. In particular, \emph{Mandersloot}, a logistics provider specializing in refrigerated transportation, provided us with technical information about its refrigerated vehicle fleet, %along with the results of the ATP-compliant technical test performed on a trailer equipped with a Thermo King A500 refrigeration unit. 
while \emph{Tcomm Telematics}, a company that specializes in IT support and sensor technologies, and \emph{ORTEC}, who designed decision support software for Mandersloot, gave us access to data on the routes traveled by Mandersloot's fleet, including driving times, sequence of customers visited, quantities and types of products unloaded along the route with associated door opening times. 

We select four routes traveled on June 1, 2023, in the Netherlands, the Czech Republic, Hungary, and Poland, respectively, as a basis for our instances. We refer to these routes as R1, R2, R3, and R4; their main characteristics are summarized in Table \ref{tab:Routes-details}. To set the external temperature along the routes, we chose Rotterdam, Prague, Budapest, and Warsaw as representative cities for the $4$ routes. We choose the $75^{\text{th}}$ percent quantiles of the long-term temperature distributions in these cities at the beginning of June to showcase the performance of the proposed policies in a slightly challenging, albeit not extreme, setting. Data on external temperatures are sourced from \url{https://weatherspark.com/}.

\begin{table}[t]
    \small
    \centering
    \begin{tabular}{lccccccc}
        \toprule
        Route & Number of & Average Driving Time & Average Stop & Route Duration & \multicolumn{3}{c}{Temperature (\unit{\kelvin})} \\
        \cline{6-8}
        & Stops & Between Stops (min) & Time (min) & (hours) & Min & Avg & Max \\
        \midrule
        R1 & 6 & 72 & 17 & 8.5 & 288 & 291 & 293 \\
        R2 & 6 & 88 & 11 & 9.8 & 288 & 293 & 296 \\
        R3 & 11 & 47 & 12 & 10.7 & 291 & 297 & 299 \\
        R4 & 6 & 75 & 12 & 8.6 & 293 & 296 & 297 \\
        \bottomrule
    \end{tabular}
    \caption{\label{tab:Routes-details} Routes description with external temperature details.}
\end{table}

We assume four different initial pallet loads and manually specify the quantities loaded and unloaded at each stop, ensuring proportionality with historical door opening times. The exact loading and unloading schedules are presented in Appendix \ref{app:handling_schedule}.

The definitions of the parameters used in our thermodynamic model are sourced from the literature and online resources. The coefficients for converting the electrical power consumed by the compressor into liters of fuel required by the unit are provided by \cite{stellingwerf2018reducing}, while \cite{kondjoyan2006review} provide a detailed analysis of heat transfer coefficients for beef under various chilling conditions. For this reason, we select meat as the type of cargo transported. A complete overview of the values and sources of these parameters can be found in Table \ref{tab:parameters} in Appendix \ref{app:case_study_parameters}.

Finally, to obtain the coefficients for the $K$ affine functions in \eqref{eq:cooling_power_approx}, we first determine the condensing pressure values of the refrigerant at each recorded ambient temperature during the routes \citep{larsen2007potential}. Based on this, we determine the values $\theta$ as described in the Appendix \ref{app:cooling_unit}. Secondly, to apply the least squares partition algorithm, we generated a uniform grid of $10^4$ data points, representing combinations of independent variables $\tcu$ and $\ttruck$. Based on these values, we then apply the methods in \cite{magnani2009convex} to approximate the function, choosing $K = 4$ partitions. We note that higher values of $K$ have no significant benefit in accuracy as measured by weighted mean absolute error (wMAPE). For example, setting $K = 9$ results in a wMAPE of $1.9$\% as compared to the $3.2$\% with $K = 4$. In the interest of computational tractability, we therefore choose the latter approximation.

Based on the parameters collected, we set up an extensive battery of $36$ test cases by varying two key parameters for each of the four routes. Firstly, we vary the refrigeration unit's capacity to determine whether a more efficient utilization of the unit could justify the purchase of a less powerful machine. The range of values for the grid search is chosen in accordance with the technical reports of Mandersloot and relevant online sources, specifically 8-10-12 [\unit{k\watt}]. Secondly, we vary the heat transfer coefficient, which is a parameter influenced by factors such as product type, shape, airflow properties (including velocity and turbulence), and packaging design. Lower values of $h$ reduce the vulnerability of the product to temperature fluctuations, leading to more stable temperatures and lower costs in both temperature deviations and energy consumption. Based on the findings proposed in \cite{kondjoyan2006review}, we select as possible values: 2-4-6 [\SI{}{W\per\meter\squared\kelvin}].

%%% SECTION 4.2. Stochastic Models
\subsection{Samples Generation}\label{ssec:samples_generation} 
In the following, we describe how we fit models for uncertain door opening times and model the uncertainty in initial product temperatures. The resulting specification of the stochastic process $\xi = (\xi_1, \dots, \xi_S)$ is subsequently used to generate scenarios that serve as inputs for the generation of the scenario lattice, which is used as a surrogate for $\xi$ in SDDP \citep{MinnerLohndorfWozabal2013, LoWo21}. A separate set of scenarios is the basis for the out-of-sample evaluation of all policies.

Since we do not have data on the initial temperatures of loaded products, we sample from a multivariate uniform distribution on $[\underline T_p, \overline T_p]$.  Note that this choice ensures that the products never arrive in a state of violation of their temperature limits. While it would be possible to relax this assumption and allow for instant violations of temperature limits after loading, this would lead to certain deviation cost in our objective, thus skewing the comparison of policies. Furthermore, it is not clear whether the transport company should be held accountable for such deviations. Additionally, we assume that the temperatures of the products loaded at the same stop are correlated with a linear correlation coefficient of $0.8$ to take into account that the loaded products were probably stored under similar storage conditions.

In order to come up with a stochastic model for door opening times, we use a data set of $11556$ records of historical stops that include information on door opening time along with the amount of product loaded (in \emph{loaded meters}) as well as the type of the cargo. We removed records where the door opening times exceeded 3 hours and retained only records with an opening time per loaded meter rate of less than 20 minutes. We also remove outliers with more than 13.2 loaded meters unloaded, which is the truck's capacity, resulting in a data set with $7312$ samples. The data set is then divided into training and test sets using stratified sampling to ensure that records belonging to the same route are not separated and that the proportion of product types remains homogeneous in both sets.

In a next step, we apply XGBoost \citep[see][]{chen2016xgboost}, a state-of-the-art machine learning technique based on gradient boosted trees, to predict door opening times $y$ using product types and handled quantities as features. To this end, we use one-hot encoded product types, loaded quantities, as well as all quadratic interactions of these basic covariates as features $X$ for the model. To avoid negative values for predicted door opening times, we estimate a model $\hat f$ for logarithmic door opening times of the form
$$ \ln y = f(X) + \varepsilon.$$
In the training phase, we use 5-fold cross-validation and a grid search to tune the hyperparameters \emph{learning rate} (best parameter: $0.17$), a \emph{max-depth} (best parameter: $4$) and \emph{number of boosting rounds} (best parameter: $19$). Considering the results expressed in minutes (not on a logarithmic scale), the resulting model has a root mean squared error of $15.5$ minutes and an $R^2$ of $0.45$. 

Finally, to sample scenarios, we use the known data on cargo type and handled quantities $x_s$ to predict a log-door opening time using our fitted model $\hat f$ and then re-sample from the empirical noise $\varepsilon$ to obtain a scenario $\exp(\hat f(x_s)+\tilde \varepsilon)$ at a specific stop. Note that the scenarios for the same stop, therefore, differ in the sampled residual $\tilde \varepsilon$.

We model door opening times and initial product temperatures as stage-wise independent, i.e., we assume that $\xi_s$ and $\xi_{s+1}$ are independent. To test this assumption for door opening times, we perform an autocorrelation analysis of the residuals $\tilde \varepsilon$. Specifically, the residuals are reordered obtaining a series of chronological sequences of residuals corresponding to the same historical routes. We than estimate a linear regression of the form 
$$ \tilde \varepsilon_t = \beta_0 + \beta_1 \tilde \varepsilon_{t-1} + \nu_t$$
where $\tilde \varepsilon_{t-1}$ is the residual preceding $\tilde \varepsilon_t$ in the tour, and we use all consecutive pairs of residuals in the data as response and covariates. The estimation produces a $99\%$ confidence interval for $\beta_1$ ranging from $-0.03$ to $0.18$, allowing us to reject the hypothesis of autocorrelation on the first lag and thereby making our assumption of stage-wise independent randomness plausible.

%%% SECTION 4.3. Numerical Results
\subsection{Numerical Results}\label{ssec:numerical_results} 
Below, we analyze the results of simulating $1000$ out-of-sample scenarios for each of the $36$ problem instances. The experiments were performed on a computer running Microsoft Windows 10 Pro with an Intel Core i7-6700K CPU operating at 4.00 GHz, featuring 4 cores and 8 logical processors, and equipped with 24 GB of RAM. The \emph{Quasar} library and the CPL solver were used to solve SDDP (SP) and RLP and perform the out-of-sample test for all policies except our heuristic strategies H1 and H2. 

We set the number of SP iterations to $450$ and create scenario lattices with $100$ nodes from stage $2$ onward based on $10^5$ scenarios of $\xi$ for lattice generation. To achieve sufficiently accurate time discretization, we set the parameter $\Delta_D = 60$ [\unit{\second}], while the number of intervals corresponding to each opening phase equals the maximum number of minutes that the doors are open in the scenarios used for lattice generation, restricting the maximum length of a time interval to one minute. 

We set $\lambda_1 = 0.5$ in the heuristic policies, which strikes a balance between lower values that might help reduce violations but lead to rather frequent switching resulting in excessive wear of mechanical components \citep{li2010optimal, HUANG2017576} and a higher value that leads to a smoother operational policy at the cost of higher violations. Furthermore, we choose $\lambda_2 = 1$ in H2, which means that we start to cool as soon as one of the products becomes too warm.

\begin{table}[t]
        \centering
        {\linespread{1}\begin{small}\selectfont
        \begin{tabular}{>{\ttfamily}lcccccccc}
            \toprule
            \textbf{} & W & h & $\Cc^{\text{CLV}}$ & $\Cc^{\text{SP}}$ &
            $\Cc^{\text{RLP}}$ & $\Cc^{\text{H2}}$  & $\Cc^{\text{H1}}$ \\
            \midrule
    \multirow{9}{*}{\begin{sideways}R1\end{sideways}}
    & \multirow{3}{*}{8} & 2 & 0.6 & 0.7 (0.2) & 3.9 (0.7) & 8.8 (1.1)  & 24.1 (2.8) \\
    & & 4 & 4.0 & 5.2 (1.0) & 46.5 (5.1)  & 91.1 (6.8) & 212.5 (14.7) \\
    & & 6 & 20.6 & 24.4 (3.6) & 148.1 (12.6)  & 214.2 (13.5) & 379.9 (22.1) \\
\cmidrule(lr){2-8}    & \multirow{3}{*}{10} & 2 & 0.6 & 0.7 (0.2) & 4.2 (0.7) & 8.6 (1.0)  & 23.3 (2.7) \\
    & & 4 & 3.9 & 5.1 (0.9) & 46.7 (5.2)  & 89.8 (6.7) & 209.5 (14.6) \\
    & & 6 & 20.2 & 24.2 (3.6) & 148.7 (12.5)  & 211.4 (13.4) & 374.9 (22.0) \\
\cmidrule(lr){2-8}    & \multirow{3}{*}{12} & 2 & 0.6 & 0.8 (0.2) & 3.7 (0.7) & 8.5 (1.0)  & 23.1 (2.7) \\
    & & 4 & 3.9 & 5.1 (0.9) & 46.5 (5.1)  & 89.6 (6.7) & 208.7 (14.6) \\
    & & 6 & 20.2 & 24.1 (3.6) & 146.6 (12.3)  & 210.9 (13.4) & 374.3 (22.0) \\
\midrule\multirow{9}{*}{\begin{sideways}R2\end{sideways}}
    & \multirow{3}{*}{8} & 2 & 0.0 & 0.1 (0.1) & 2.4 (0.4) & 3.4 (0.4)  & 11.1 (1.4) \\
    & & 4 & 0.2 & 1.2 (0.3) & 34.4 (4.2)  & 36.7 (4.0) & 95.5 (9.8) \\
    & & 6 & 3.7 & 8.1 (1.8) & 88.6 (9.3)  & 112.9 (10.2) & 199.3 (17.0) \\
\cmidrule(lr){2-8}   & \multirow{3}{*}{10} & 2 & 0.0 & 0.1 (0.1) & 2.1 (0.3) & 3.2 (0.4)  & 10.3 (1.4) \\
    & & 4 & 0.2 & 1.1 (0.3) & 32.5 (4.2)  & 35.4 (3.9) & 92.8 (9.6) \\
    & & 6 & 3.6 & 8.1 (1.8) & 87.7 (9.4)  & 110.8 (10.1) & 197.3 (16.9) \\
\cmidrule(lr){2-8}   & \multirow{3}{*}{12} & 2 & 0.0 & 0.2 (0.1) & 2.1 (0.4) & 3.1 (0.4)  & 10.0 (1.3) \\
    & & 4 & 0.2 & 1.1 (0.3) & 32.8 (4.1)  & 35.0 (3.9) & 91.5 (9.5) \\
    & & 6 & 3.6 & 10.7 (2.3) & 89.0 (9.2)  & 109.8 (10.0) & 195.7 (16.9) \\
\midrule\multirow{9}{*}{\begin{sideways}R3\end{sideways}}
    & \multirow{3}{*}{8} & 2 & 29.1 & 32.1 (2.8) & 48.6 (4.2) & 273.8 (11.9)  & 759.8 (21.4) \\
    & & 4 & 225.9 & 249.8 (16.4) & 341.9 (20.8)  & 1151.3 (36.9) & 2191.4 (49.6) \\
    & & 6 & 745.7 & 780.6 (38.3) & 818.4 (39.7)  & 2158.0 (61.4) & 3051.7 (66.4) \\
\cmidrule(lr){2-8}   & \multirow{3}{*}{10} & 2 & 24.7 & 27.6 (2.5) & 41.7 (3.7) & 243.9 (11.1)  & 691.1 (2
1.0) \\
    & & 4 & 207.0 & 230.0 (15.6) & 304.7 (19.4)  & 1097.6 (36.0) & 2167.9 (48.6) \\
    & & 6 & 677.6 & 711.6 (35.8) & 752.6 (37.3)  & 1974.2 (57.9) & 3037.8 (66.0) \\
\cmidrule(lr){2-8}   & \multirow{3}{*}{12} & 2 & 22.5 & 25.3 (2.3) & 39.2 (3.7) & 229.9 (10.7)  & 657.9 (2
0.5) \\
    & & 4 & 198.0 & 220.9 (15.2) & 319.0 (20.0)  & 1068.7 (35.5) & 2122.3 (48.3) \\
    & & 6 & 660.4 & 693.9 (35.2) & 737.9 (37.0)  & 1941.0 (57.4) & 3001.6 (65.7) \\
\midrule\multirow{9}{*}{\begin{sideways}R4\end{sideways}}
    & \multirow{3}{*}{8} & 2 & 1.3 & 1.3 (0.2) & 5.9 (0.8) & 9.0 (0.9)  & 23.4 (2.4) \\
    & & 4 & 6.3 & 7.4 (1.0) & 61.3 (5.1)  & 80.1 (5.7) & 192.0 (12.6) \\
    & & 6 & 23.2 & 29.7 (3.5) & 136.6 (9.8)  & 231.4 (12.5) & 400.8 (19.9) \\
\cmidrule(lr){2-8}    & \multirow{3}{*}{10} & 2 & 1.2 & 1.2 (0.2) & 6.2 (0.9) & 8.1 (0.8)  & 21.2 (2.3) \\
    & & 4 & 6.0 & 7.1 (1.0) & 54.8 (4.8)  & 73.5 (5.4) & 178.4 (12.4) \\
    & & 6 & 22.6 & 28.8 (3.4) & 130.7 (9.5)  & 221.1 (12.2) & 388.3 (19.7) \\
\cmidrule(lr){2-8}   & \multirow{3}{*}{12} & 2 & 1.1 & 1.2 (0.2) & 5.2 (0.8) & 7.8 (0.8)  & 20.3 (2.2) \\
    & & 4 & 5.9 & 7.0 (0.9) & 55.3 (4.8)  & 72.0 (5.4) & 174.6 (12.2) \\
    & & 6 & 22.3 & 28.7 (3.4) & 139.3 (10.0)  & 217.8 (12.1) & 382.8 (19.5) \\
\midrule
        \end{tabular}
        \end{small}}
        \caption{\label{tab:results} Average cost [\unit{\kelvin\cdot\min}](standard errors in brackets) of $1000$ out of sample tests for varying routes, cooling capacity [k\unit{\watt}], and heat transfer coefficient $h$  [\SI{}{W\per\meter\squared\kelvin}].}
\end{table}

Table \ref{tab:results} reports the results on the out-of-sample temperature violations for all $36$ problem instances and policies. The clairvoyant policy (CLV) serves as a perfect-foresight benchmark which is computed by the average cost of the deterministic problems computed for every scenario, thus assuming full knowledge of the realizations of $\xi$. As such, CLV is not implementable and provides an optimistic lower bound on what can be achieved by any implementable policy in our setting.

%%% SECTION 4.3.1 Overview of Results
\subsubsection{Overview of Results}
Inspecting Table \ref{tab:results}, we see that the out-of-sample costs of the policies are consistently ordered according to $\Cc^{\text{H1}}>\Cc^{\text{H2}}>\Cc^{\text{RLP}} > \Cc^{\text{SP}} > \Cc^{\text{CV}}$. This reflects the varying degree of sophistication as well as the amount of information used by the policies. In particular, we observe that SP consistently and significantly outperforms the best available deterministic benchmark RLP, leading to a sizable value of the stochastic solution in most instances. Furthermore, the differences between the lower bound CLV and SP are small in most instances, indicating that the stochastic solution performs very well. Inspecting the solutions for H1 and H2, we see that these simpler heuristic policies usually perform much worse. The maximum training time is $67$ minutes for the SP policy and for a large majority of $92\%$ of cases it is less than an hour, making it a feasible solution for real-world scenarios.

\begin{figure}[t]
    \centering
    \includegraphics[width=1\textwidth]{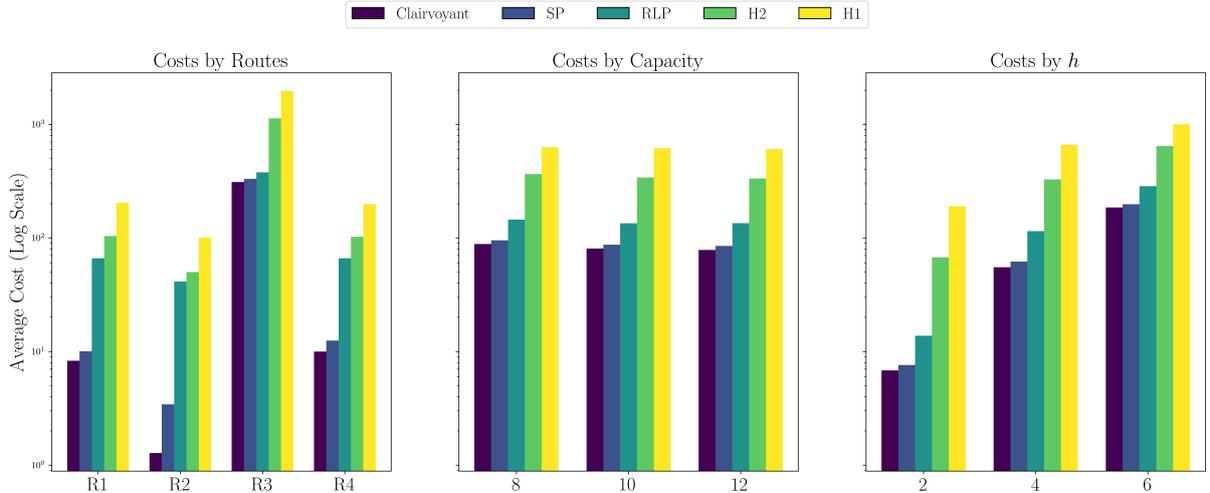}
    \caption{Logarithmic scale representation of the average performance of each policy [\unit{\kelvin\cdot\min}], varying route (left), refrigeration unit capacity [k\unit{\watt}
] (center), and heat transfer coefficient [\SI{}{W\per\meter\squared\kelvin}] (right).}
    \label{fig:parameters_results}
\end{figure}

Figure \ref{fig:parameters_results} averages the results and groups them according to each of our parameters. From the graph on the left, we observe that route R3 is the most challenging, as it has the highest cost, but also (relatively) the smallest difference between the policies. This stems from the fact that this route has the highest outdoor temperatures, as well as the highest number of stops, which are additionally separated by relatively short driving intervals. These factors, in combination, make it difficult to manage the temperature of loaded products such that even CLV struggles to achieve low costs, thus narrowing the potential for effective cooling. However, substantial cost reductions can still be achieved compared to the heuristic rules. 

An inspection of the second panel of Figure \ref{fig:parameters_results} reveals that surprisingly, the capacity of the refrigeration unit has little impact on deviations. This suggests that the efficient use of the 8 kW unit yields similar performance to the $12$ kW unit, potentially making the latter's purchase unnecessary. In fact, looking at the detailed results in Table \ref{tab:results}, we see that the cooling capacity only has a significant effect for R3, that is, for the route characterized by short driving intervals and many stops, although even in this case the effect is not significant.

The last panel in Figure \ref{fig:parameters_results} shows the impact of the coefficient $h$. Clearly, this parameter, which represents the heat conductivity of the loaded products, has a significant impact, indicating that investments in advanced packaging pay off. Looking at the difference between policies, we see that SP and CLV perform better relative to less sophisticated policies for lower values of $h$. This is because for higher values of $h$, products warm up quicker, and therefore, there are more \emph{forced} violations due to long door opening times. In contrast, for small $h$, sophisticated policies manage to avoid most of the deviations and, therefore, more clearly outperform simpler heuristics.

%%% SECTION 4.3.2 Optimal Decisions
\subsubsection{Optimal Decisions}
\begin{table}[t]
    \centering
    \small
    \renewcommand{\arraystretch}{1.2}
    \setlength{\tabcolsep}{6pt}
    \begin{tabular}{lcccccccccccc}
        \toprule
        \multirow{2}{*}{} & \multicolumn{3}{c}{SP} & \multicolumn{3}{c}{RLP} & \multicolumn{3}{c}{H2} &
 \multicolumn{3}{c}{H1} \\
        \cmidrule(lr){2-4} \cmidrule(lr){5-7} \cmidrule(lr){8-10} \cmidrule(lr){11-13}
                              & $\overline{A}$ & $\overline{W}$ & $\overline{VS}$
                              & $\overline{A}$ & $\overline{W}$ & $\overline{VS}$
                              & $\overline{A}$ & $\overline{W}$ & $\overline{VS}$
                              & $\overline{A}$ & $\overline{W}$ & $\overline{VS}$ \\
        \midrule
        R1 & 0.85 & 3834 & 30 & 0.66 & 2827 & 182 & 0.33 & 6381 & 282 & 0.31 & 6479 & 295 \\
        R2 & 0.84 & 3410 & 11 & 0.79 & 2418 & 123 & 0.33 & 7068 & 145 & 0.32 & 7129 & 157 \\
        R3 & 0.95 & 6219 & 343 & 0.97 & 5710 & 350 & 0.53 & 7910 & 691 & 0.41 & 8412 & 806 \\
        R4 & 0.88 & 4536 & 47 & 0.87 & 3528 & 238 & 0.34 & 7906 & 313 & 0.32 & 8075 & 331 \\
        \bottomrule
    \end{tabular}
    \caption{\label{tab:decisions} Cooling decisions statistics obtained from out-of-sample simulations performed on 36 different instances. For each policy and each route, the average percentage of time the refrigeration unit is active ($\overline{A}$), the average power consumed ($\overline{W}$),    and the average number of scenarios in which violations occurred ($\overline{VS}$) are reported.}
\end{table}

Next, we turn our attention to the cooling decisions. We start by discussing some operational metrics averaged over all $h$ and the refrigeration unit capacities as reported in Table \ref{tab:decisions}. The results clearly show that for SP and RLP, the refrigeration unit is consistently switched on longer than for H1 and H2, but at the same time, significantly less average power is used. In particular, for SP, the refrigeration unit is on average on for more than 80\% of the time, but at a relatively low setting, ensuring a continuous but less aggressive cooling strategy than the other policies. As was expected, the number of violations for the two heuristics is significantly higher than those for RLP and especially SP. 

\begin{figure}[t]
    \centering    \includegraphics[width=0.9\textwidth]{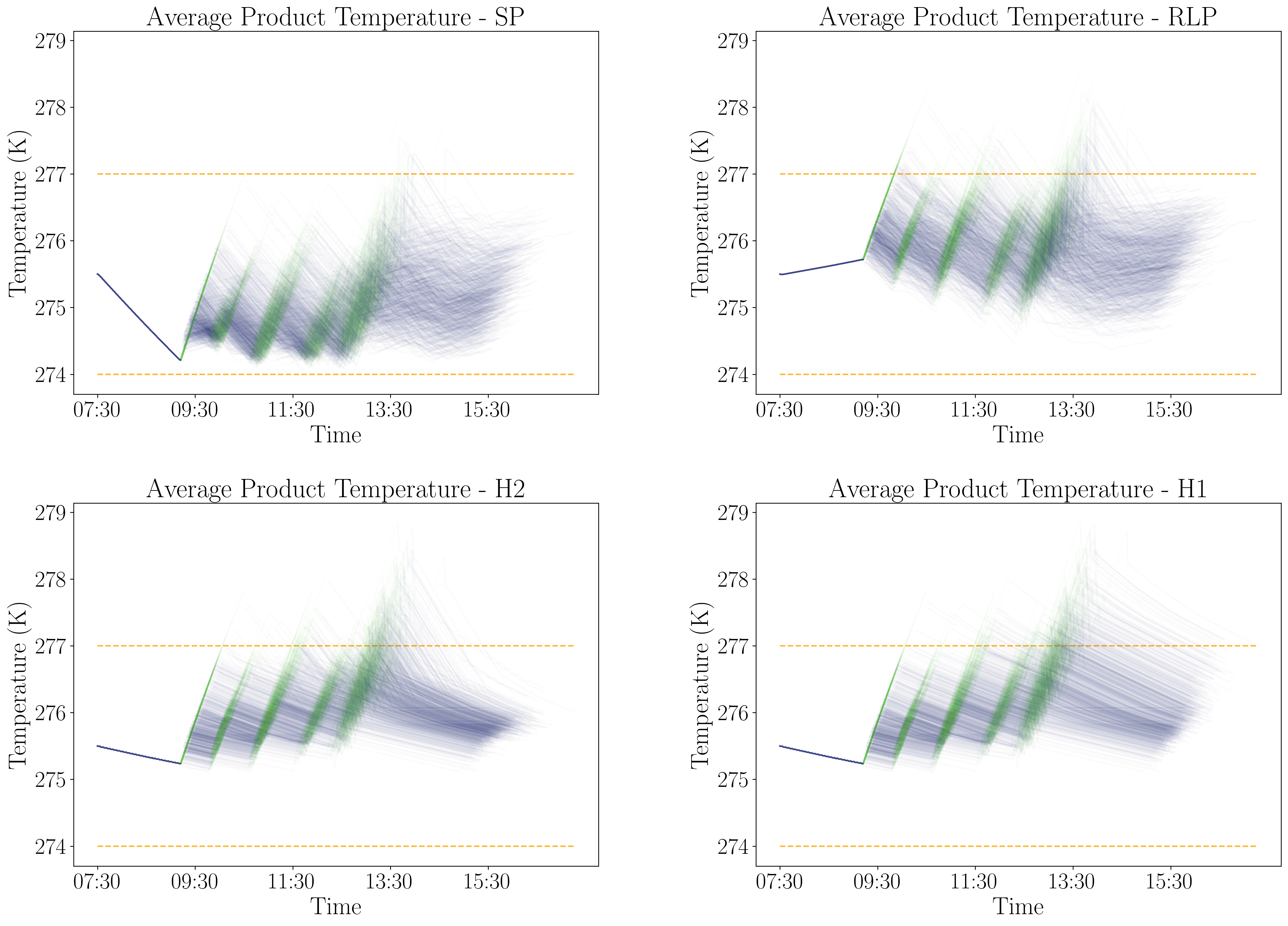}
    \caption{Evolution of average product temperatures for all out-of-sample scenarios for route R1 with unit capacity $12$ and $h=4$. The $x$ axis represents actual clock time, which results in different trip lengths due to the random door opening duration. Green periods correspond to handling operations with open doors, while blue periods represent driving.}
    \label{fig:POJ-03_results}
\end{figure}

Finally, Figure \ref{fig:POJ-03_results} contains detailed graphs of temperatures and cooling decisions for R1 that compare the evolutions of the average product temperatures along the routes for all $1000$ scenarios and $4$ policies. Looking at the graphs, the policies mainly seem to differ in how early they start cooling. While H1 and H2 are myopic and therefore cool rather late and independent of stops, RLP and SP anticipate stops and engage in preemptive cooling, building up a thermal buffer in preparation for anticipated warming during door opening. 

When comparing RLP and SP, it is evident that SP cools much more aggressively before a stop, i.e., acts more conservatively. This clearly is the case since RLP plans for the expected scenario, while SP sees all possible stop durations, including longer ones, which, because of the asymmetric cost structure, leads to more cooling. This becomes evident when looking at the first stop: while RLP does not cool at all, since the expected door opening time would not lead to any violations of temperature bounds, SP also takes into consideration extreme scenarios in the first stop as well as in the stops after that and decides to cool down the products quite substantially in the beginning, thereby building up a thermal buffer that makes it possible to complete the trip with nearly no violations in all scenarios. Another advantage of precooling is that the refrigeration unit is more effective at the beginning of the trip since the outside temperatures are lower. Hence, we identify intelligent anticipative precooling as the main feature of a successful cooling policy.

%%% SECTION 4.3.3 The Value of Information
\subsubsection{The Value of Information}
\begin{figure}[t]
    \centering
    \includegraphics[width=1\textwidth]{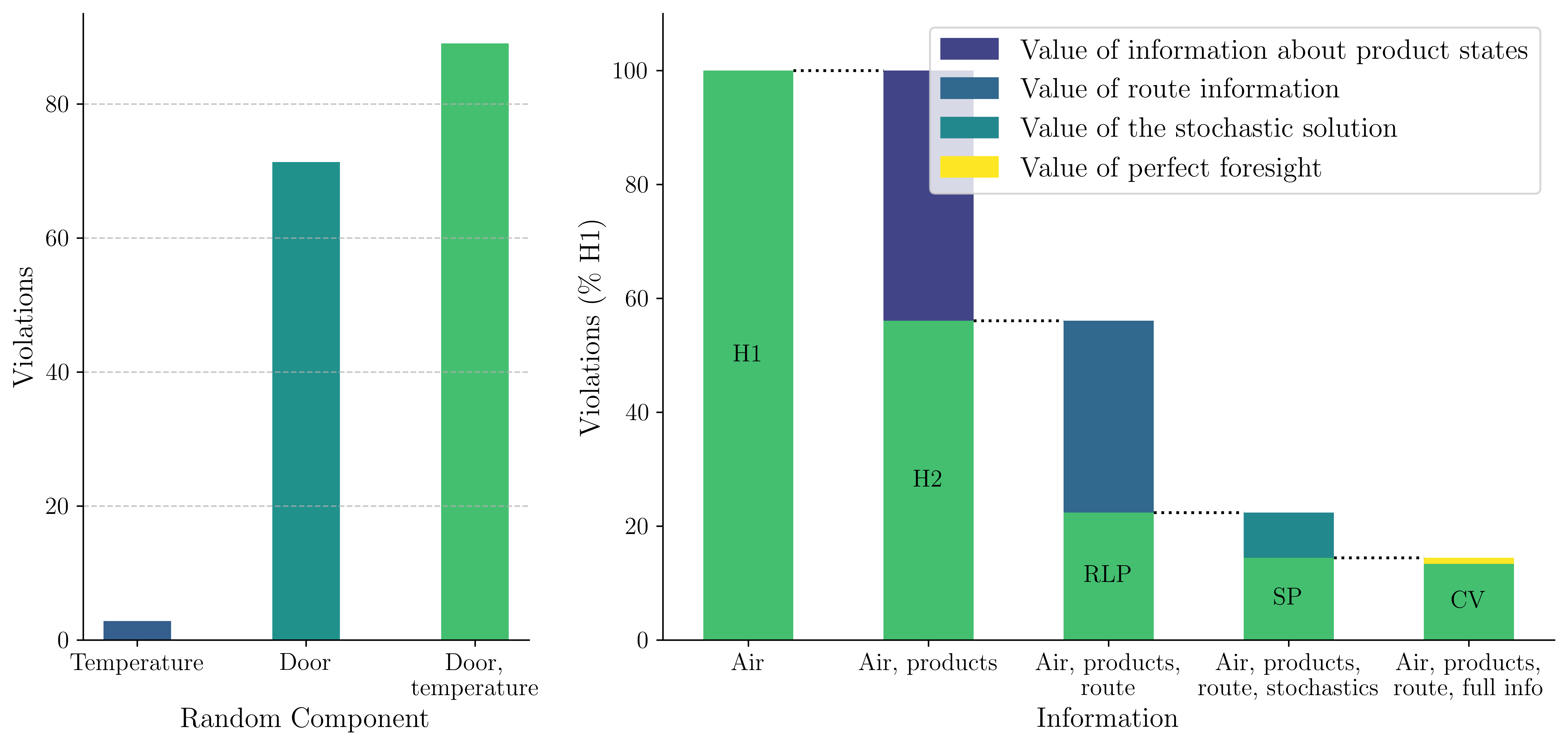}
    \caption{Effect of random components on average violations [\unit{\kelvin\cdot\min}] obtained by SP (left), value added associated with the information available to each policy, expressed as a percentage of average H1 violations (right).}
    \label{fig:results_routes}
\end{figure}
In this section, we investigate the influence of information and randomness on the solution of the problem. We start by discussing how the two random factors, door opening times and initial temperatures, influence expected cost. In Figure \ref{fig:results_routes} on the left, we compare the average cost of SP in a problem with only random opening times and deterministic initial temperatures, as well as a problem with random initial temperatures but deterministic opening times, with the model with the full randomness. For the model with partial randomness, we set the respective other variable to its expectation. The plot reveals that the door opening times induce a much higher cost than the initial temperatures in case these factors are used as the only random factors. Additionally, we observe that the effect of the two sources of randomness is not additive, i.e., the model where both factors are random has a higher cost than the two reduced models combined. Hence, the effects of the two random factors seem to mutually reinforce each other.

In the right part of Figure \ref{fig:results_routes}, we study the average costs of our four policies aggregated over all cases and interpret the differences as the value of different types of information. The difference between H1 and H2 can be interpreted as the value of having information about the products' temperatures versus only measuring air temperatures. The next step to RLP signifies the value of knowing the temporal structure of the route as well as the thermodynamic interactions in the trailer. The step from RLP to SP is the value of the stochastic solution, while the difference between SP and CLV is the value of perfect foresight. When looking at the magnitudes of the jumps, a decreasing marginal value of information becomes apparent. In particular, the greatest impact comes with measuring the products' thermal states, which is an indication that investment in advanced measurement devices can lead to large cost savings. Furthermore, the relatively large jump from H2 to RLP indicates that knowing the route of the truck represents a major advantage that should be used by any cooling policy, as it allows preemptive cooling. Lastly, the difference between SP and CLV is surprisingly small, indicating that SP is close to optimal. 

%%% SECTION 4.3.4 The Quality-Fuel Consumption Trade Off
\subsubsection{The Quality-Fuel Consumption Trade-Off}
\begin{figure}[t]
    \centering    \includegraphics[width=1\textwidth]{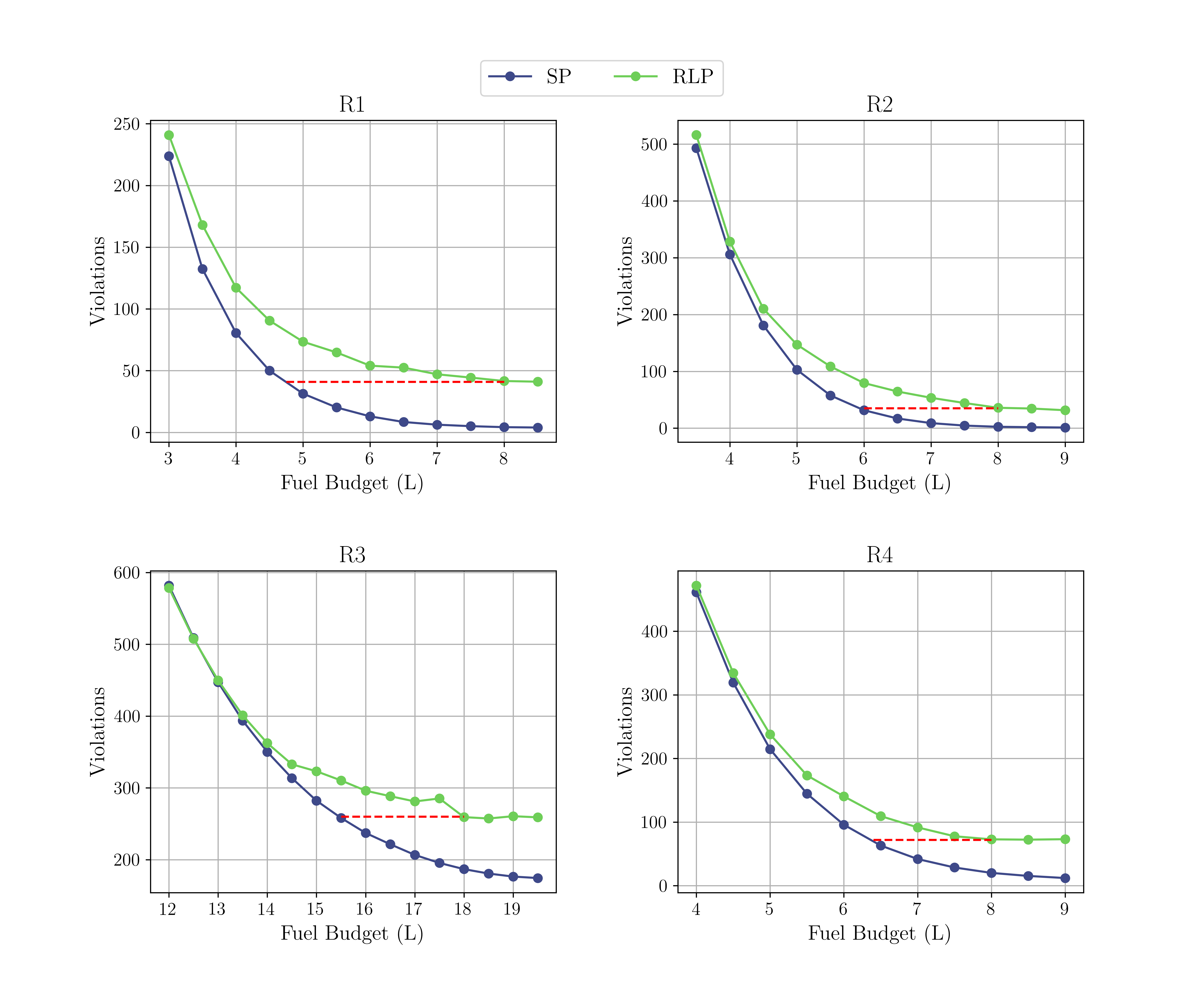}
    \caption{Comparison of results of model \eqref{eq:fuel_budget_constraint} with fuel budget constraint. For each route, average violations [\unit{\kelvin\cdot\min}] are shown as the liters of fuel available to the unit change.}
    \label{fig:results_fuelBudget}
\end{figure}

In this section, we discuss the model variant in which we induce a limited fuel budget by adding \eqref{eq:fuel_budget_constraint} to the model. In particular, we are interested in how the RLP performs compared to SP when we impose an upper limit on fuel consumption. We do so in Figure \ref{fig:results_fuelBudget} which shows efficient frontiers depicting the trade-off between the fuel budget and the violations for the two policies. For the purpose of this analysis, we fix $h=4$, the capacity of the cooling unit to $12$ and set the penalty $\mu$ to $10000$, which is large enough to result in policies without any violation of the fuel budget.

Looking at Figure \ref{fig:results_fuelBudget}, we notice that, consistent with previous results, the efficient frontiers for SP are below those for RLP. This effect becomes more pronounced with increasing fuel budget, since for a lower budget the amount of unavoidable violations, which are independent of the chosen policy, grows. As the fuel budget increases, the policies have enough fuel to avoid more violations, and the different efficiency of RLP and SP in doing so starts to show.

This comparison reveals that the potential for fuel savings is significant. Depending on the acceptable level of violations and the route, the difference in fuel consumption is as large as $40\%$ as exemplified by the dotted red lines. Looking, for example, at the first panel, we see that if violations of $40$ are acceptable to the decision maker, SP can achieve this goal with $4.8$ liters of fuel, while SP requires $8$. This analysis highlights the potential benefits of a sophisticated cooling policy on the energy efficiency of the industry. Therefore, policies like SP show a large potential to instantly make the industry more sustainable without relying on new technologies or large investments in updated equipment and, therefore, can significantly contribute to reaching the climate targets of the transport industry.

%%%%%%%%%%%%%%%%%%%%%
%%% 5. Conclusion %%%
%%%%%%%%%%%%%%%%%%%%%
\section{Conclusion} \label{sec:conclusion}
We introduce a novel model to optimize the working regime of the refrigeration unit installed in cooling trucks to improve the reliability and efficiency of road transportation for perishable products. Our model incorporates thermodynamic interactions between the refrigeration unit, the internal and external environment, and the products while accounting for the main stochastic factors. Our results demonstrate that our model can be solved optimally in less than an hour, while outperforming all the benchmarks: two myopic heuristics commonly used in practice (H1, H2) and a rolling lookahead policy (RLP).

By analyzing the structure of our policy's optimal solutions, we derive best practices that can be easily implemented in everyday operations. In particular, our findings suggest that decision-makers should apply preemptive cooling before door openings and measure product temperatures, not just air temperatures, to make informed decisions.

The utility of this insight can be demonstrated by optimizing the parameter $\lambda_2$ of $H_2$ for precooling. For this purpose, we tested H2 on each route for every $\lambda_2$ value between $0.1$ and $1$, in increments of $0.1$, to determine the best setting for every route based only on the scenarios used to build the scenario lattice. The optimal values for each route were found to be $0.5, 0.7, 0.1$ and $0.4$, respectively. Since on average these $\lambda_2$s are substantially lower than $1$, this indicates that it is beneficial to keep a thermal buffer based on the product temperatures and to start cooling before the products approach their temperature limit. Using the optimized parameters results in an average cost reduction of $67$\% compared to the case with $\lambda_2 = 1$. Furthermore, the optimized H2 outperforms the RLP strategy by $17$\% on average, indicating that heuristic policies that are designed to structurally mimic the optimal solution could come close to SP at a much lower computational cost and algorithmic complexity. 

Future research may focus on improving the thermodynamic model by including additional aspects such as solar radiation and humidity and designing new heuristic policies that better make use of our findings. Furthermore, it would be interesting to consider a more integrated problem in which the refrigeration policy is coupled with inventory and vehicle routing and to investigate how these decisions impact each other.

\bibliographystyle{abbrvnat}
\bibliography{references.bib}

\newpage

\appendix

\setcounter{table}{0}
\setcounter{section}{0}
\setcounter{equation}{0}
\setcounter{figure}{0}
\renewcommand{\thetable}{A\arabic{table}}
\renewcommand{\theequation}{A\arabic{equation}}
\renewcommand{\thesubsection}{A.\arabic{subsection}}
\renewcommand{\thesection}{A.\arabic{section}}
\renewcommand{\thefigure}{A\arabic{figure}}

\section{Cooling Unit} \label{app:cooling_unit}
As mentioned in \cite{larsen2007potential}, parameters of equation \eqref{eq:cooling_power} depend directly on the refrigerant used (via $a$ and $b$ coefficients), the properties of the evaporator (via $\gamma$), and the external temperature (via $P_{ci}$):
\begin{align}
    \theta_{i1} &= -\frac{\gamma a_{21}}{a_1 P_{ci} + b_1}, \\
    \theta_{i2} &= -\frac{\gamma(a_{22} P_{ci} + b_2)}{a_1 P_{ci} + b_1}, \\
    \theta_{i3} &= \frac{\gamma a_{21}}{a_1 P_{ci} + b_1}, \\
    \theta_{i4} &= \frac{\gamma(a_{22} P_{ci} + b_2)}{a_1 P_{ci} + b_1}.
\end{align}

Where $P_{ci}$ is the fluid pressure at the condenser which depends on the temperature of the external environment $T^{\text{ext}}_i$ and therefore varies with time.%of the refrigeration system, which is controlled according to \(P_c = P_\text{sat}(T_c)\), \(P_c\) becomes a function of \(T_a\) ($T_c = T_a + cost$) directly.

The following coefficient values are used in the approximations of the refrigerant properties for R134a \citep{larsen2007potential}
\begin{align}
    a_1 = -6.7014 \cdot 10^3, \; b_1 = 2.1729 \cdot 10^5, \; a_{21} = -9.0496 \cdot 10^2, \; a_{22} = 3.367 \cdot 10^3, \; b_2 = -3.152 \cdot 10^3.
\end{align}

\subsection{Case Study Parameters} \label{app:case_study_parameters}
\begin{table}[ht!]
    \centering
    \begin{tabular}{l l l l l}
        \toprule
        \textbf{Parameter} & \textbf{Value} & \textbf{Unit} & \textbf{Description} & \textbf{Source} \\ 
        \midrule
        $c'$ & 1005.6 & \SI{}{J\per\kilogram\kelvin}   & Internal air specific heat & \cite{CuriosityFluids} \\ 
        $v$ & 2.6 & \unit{\kilogram\per\second}  & Air infiltration rate & \cite{zhang2018computational} \\ 
        $A$ & $153.9$ & \unit{\meter\squared} & Trailer surface & Mandersloot \\ 
        $A'$ & $146.9$ & \unit{\meter\squared} & Walls surface & Mandersloot\\ 
        $A_p$ & $1.9$ & \unit{\meter\squared} & Product surface & Europallet\\ 
        U & 0.4 & \SI{}{W\per\meter\squared\kelvin} &   Wall’s transmittance & Mandersloot \\  
        $\gamma$ & $300$ & \SI{}{W\per\meter\squared\kelvin} & Evaporator's transmittance  & \cite{TheEngineeringToolBox} \\ 
        $ \overline{W}$ & 8-10-12 & \unit{k\watt} & Unit capacity & Mandersloot\\
        $ \underline{T_p},\ \overline{T_p}$ & 274-277 & \unit{\kelvin} & Temperature range & \cite{AUNG2014198} \\
        $h_p$ & 2-4-6 & \SI{}{W\per\meter\squared\kelvin}  & Product heat transfer coefficient & \cite{kondjoyan2006review} \\ 
        $\sigma$ & $9.504*10^3$ & \unit{k\joule\per\liter}  & Fuel conversion factor & \cite{stellingwerf2018reducing} \\   
        $C_p$ & $7.8*10^2$ & \unit{k\joule\per\kelvin}  & Product heat capacity & \cite{TheEngineeringToolBox} \\  
        $C$ & $120.513$ & \unit{k\joule\per\kelvin}  & Internal air heat capacity & Mandersloot \\   
        \bottomrule
    \end{tabular}
    \caption{\label{tab:parameters} Summary of the parameters employed in our model. }
\end{table}

\subsection{Handling Schedule} \label{app:handling_schedule}
\begin{table}[ht!]
\centering
\begin{tabular}{ccccccccc}
    \toprule
    \multirow{2}{*}{Stop} & \multicolumn{2}{c}{R1} & \multicolumn{2}{c}{R2} & \multicolumn{2}{c}{R3} & \multicolumn{2}{c}{R4} \\
    \cmidrule(lr){2-3}
     \cmidrule(lr){4-5}
     \cmidrule(lr){6-7}
     \cmidrule(lr){8-9}
         & In  & Out & In  & Out & In  & Out & In  & Out \\
    \midrule
    Depot  & 25  & 0   & 20  & 0   & 18  & 0   & 15  & 0   \\
    1      & 2   & 8   & 2   & 1   & 2   & 1   & 2   & 1   \\
    2      & 2   & 0   & 2   & 1   & 2   & 2   & 2   & 4   \\
    3      & 2   & 3   & 2   & 4   & 2   & 0   & 2   & 1   \\
    4      & 2   & 1   & 2   & 3   & 2   & 1   & 2   & 3   \\
    5      & 2   & 8   & 2   & 1   & 2   & 2   & 2   & 3   \\
    6      & 0   & 15  & 0   & 20  & 2   & 3   & 0   & 13  \\
    7      & -   & -   & -   & -   & 2   & 7   & -   & -   \\
    8      & -   & -   & -   & -   & 2   & 2   & -   & -   \\
    9      & -   & -   & -   & -   & 2   & 3   & -   & -   \\
    10     & -   & -   & -   & -   & 2   & 2   & -   & -   \\
    11     & -   & -   & -   & -   & 0   & 15  & -   & -   \\
    \bottomrule
\end{tabular}
\caption{\label{tab:load-unload} Pallets loaded and unloaded at each stop for all routes.}
\end{table}

\begin{comment}
\begin{table}[t]
%\scriptsize
\centering
%\caption{Technical parameters}
\caption{To be removed}
\label{tab:parameters}
\begin{tabular}{>{\raggedright\arraybackslash}p{10cm} >{\centering\arraybackslash}p{2cm} >{\raggedright\arraybackslash}p{3cm}}
\toprule
\textbf{Parameter} & \textbf{Value} & \textbf{Unit of measurement} \\ 
\midrule
Specific heat of air at constant pressure in the compartment & 1005.6 & J/(kg\ K) \\ 
Specific heat of air at constant pressure out the compartment & 1006.1 & J/(kg\ K) \\ 
Density of air in the compartment & 1.29 & kg/m\textasciicircum{}3 \\ 
Density of air out the compartment & 1.20 & kg/m\textasciicircum{}3 \\ 
Compartment length & 13.2 & m \\ 
Compartment height & 2.7 & m \\ 
Compartment width & 2.6 & m \\ 
Global heat transfer coefficient of the thermal insulation panel material & 0.4 & W/(m\textasciicircum{}2\ K) \\ 
Energy content of fuel & 8.8 & kWh/L \\ 
Conversion efficiency of chemical to refrigeration energy & 0.3 & kW/L \\  
Pallet length & 0.8 & m \\ 
Pallet height & 2 & m \\ 
Pallet width & 1.2 & m \\ 
Pallet mass & 300 & kg \\ 
Meat specific heat & 2600 & J/(kg\ K) \\  
Duration time interval $\Delta$ & 60 & s \\ 
\bottomrule
\end{tabular}
\end{table}
\end{comment}

\end{document}

%% file: new_command.tex
\newcommand{\mcl}[1]{\mathcal{#1}}

\newcommand{\tprod}{T} %products temperature
\newcommand{\ttruck}{T^{\text{air}}} %temperature of the air inside the truck
\newcommand{\tamb}{T^{\text{ext}}} %external temperature (ambient temperature)
\newcommand{\tcu}{T^{\text{cu}}} %Te in the paper, evaporation temperature